# TOPOLOGICAL JACOBI FORMS

TILMAN BAUER AND LENNART MEIER

ABSTRACT. As a generalization of the ring spectrum of topological modular forms, we construct a graded ring spectrum of topological Jacobi forms, $\mathrm{TJF}_*$. This is constructed as the global sections of a sheaf of $E_\infty$-ring spectra on the stacky universal elliptic curve using circle-equivariant TMF. Complete calculations of its homotopy at odd primes and partial results at $p = 2$ are given.

## Contents

Modular forms came to prominence in topology through the elliptic genera defined by Ochanine and Witten [Och09, HBJ92]. The Ochanine genus was later seen as a specialization of the 2-variable elliptic genus on almost-complex manifolds and studied extensively by algebraic geometers [BL00, BL05]. On manifolds with vanishing first Chern class, it takes values in *Jacobi forms*. As introduced in [EZ85], these are complex functions of two variables, subject to certain transformation conditions that depend on an index and a weight.

Lifting the Ochanine and Witten genera to maps of spectra was a major motivation for defining elliptic cohomology and topological modular forms TMF. The resulting map $\mathrm{MString} \to \mathrm{TMF}$ [AHR10] has led to new insights in manifold topology [Kra21, BS24].

To replicate a similar success story for the 2-variable elliptic genus, we define spectra $\mathrm{TJF}_m$ of *topological Jacobi forms* for every index $m$. Since the first versions of our paper have circulated, Lin–Yamashita [LY24] have shown that the 2-variable elliptic genus indeed lifts to maps of spectra, taking the form $\mathrm{MTSU}(m) \to \mathrm{TJF}_m$, and this has led to new divisibility result on Euler characteristics.

---

Upon complexification, $\pi_*\mathrm{TJF}_m$ agrees with the vector space of Jacobi forms of index $\frac{m}{2}$. But like $\pi_*\mathrm{TMF}$, the groups $\pi_*\mathrm{TJF}_m$ contain interesting 2- and 3-torsion. We give a complete calculation for $(\pi_*\mathrm{TJF}_m)_{(3)}$ and partial information about $(\pi_*\mathrm{TJF}_m)_{(2)}$.

Besides its role as a recipient of the 2-variable elliptic genus, TJF is interesting due to its interpretation as $\mathrm{RO}(\mathbf{S}^1)$-graded or twisted $\mathbf{S}^1$-equivariant TMF. Via the Segal–Stolz–Teichner program, it conjecturally thus classifies certain forms of anomalous supersymmetric quantum field theories with $U(1)$-symmetry [BE24]. On the purely mathematical side, the connection to equivariant TMF allows us to express $\mathrm{TJF}_m$ in terms of stunted projective spaces. Forthcoming work of Lin–Tominaga–Yamashita uses this to calculate $C_n$-equivariant TMF for $n = 2, 3$.

A more precise statement of our results and definitions is given in the following section.

*Acknowledgments.* We thank Mayuko Yamashita for her interest, Laurent Smits for his comments on a preliminary version of this manuscript and Sil Linskens for suggesting the ring structure on $\mathrm{TJF}_0$.

We also thank the Institute Mittag-Leffler for its hospitality during the initial phase of this project in Spring 2022. The second author would like to thank the Isaac Newton Institute for Mathematical Sciences, Cambridge, for support and hospitality during the programme Equivariant Homotopy theory in Context where work on this paper was undertaken. This work was supported by EPSRC grant no EP/R014604/1. The second author was supported by the NWO grant VI.Vidi.193.111a.

1. Definitions and main results

Let $\mathcal{H} = \{\tau \in \mathbf{C} \mid \mathrm{Im}(z) > 0\}$ denote the upper half plane. As is usual in the modular forms literature, we write

$$e(x) = e^{2\pi i x} \text{ for } x \in \mathbf{C}, \quad q = e(\tau) \text{ for } \tau \in \mathcal{H} \quad \text{and} \quad \zeta^r = e(rz) \text{ for } z \in \mathbf{C}, r \in \mathbf{R}.$$

Complex-valued Jacobi forms were first defined and studied in some depth in [EZ85] as a common generalization of modular forms and elliptic functions such as theta series. We follow [GW20] in their naming convention for the various versions of Jacobi forms.

**Definition** ([EZ85, CD19]). A *weakly holomorphic Jacobi form of weight $k \in \mathbf{Z}$ and index $m \in \frac{1}{2} \cdot \mathbf{Z}$* is a holomorphic function

$$\phi \colon \mathcal{H} \times \mathbf{C} \to \mathbf{C}$$

satisfying the following transformation properties for $\begin{pmatrix} a & b \\ c & d \end{pmatrix} \in \mathrm{SL}_2(\mathbf{Z})$ and $(\lambda, \mu) \in \mathbf{Z}^2$:

(1.1) $$\phi\Big(\frac{a\tau + b}{c\tau + d}, \frac{z}{c\tau + d}\Big) = (c\tau + d)^k e\Big(\frac{mcz^2}{c\tau + d}\Big)\phi(\tau, z)$$

and

(1.2) $$\phi(\tau, z + \lambda\tau + \mu) = e(m(-\lambda^2\tau - 2\lambda z + \lambda + \mu))\phi(\tau, z)$$

and such that there exists a Fourier expansion

(1.3) $$\phi(\tau, z) = \sum_{n \geq -N} \sum_{r \in \mathbf{Z} + m} c(n, r) q^n \zeta^r \quad \text{for some } c(n, r) \in \mathbf{C}.$$

A weakly holomorphic Jacobi form is called a *weak Jacobi form* if $c(n,r) = 0$ for $n < 0$, and a *(holomorphic) Jacobi form* if additionally, the stronger cusp holomorphicity condition

$$c(n,r) = 0 \quad \text{if } r^2 > 4mn \tag{1.4}$$

holds. We use a nonstandard bigrading that is more suited to our homotopy theoretic applications and call $n = 2k + 4m$ the *dimension* of $\phi$ and set $|\phi| = (n, 2m)$.

Denote the vector space of weakly holomorphic Jacobi forms of index $m$ and dimension $n$ by $\mathrm{JF}^{\mathbf{C}}_{n,2m}$ and the subspaces of weak and of holomorphic Jacobi forms by $\mathrm{jF}^{\mathbf{C}}_{n,2m}$ and $\mathrm{jf}^{\mathbf{C}}_{n,2m}$, respectively.

All versions of Jacobi forms are bigraded $\mathbf{C}$-algebras with respect to the pointwise addition and multiplication of functions.

In analogy with modular forms, a $\mathbf{Z}$-valued version of weak Jacobi forms can be defined using the language of moduli stacks. Let $\mathcal{M}_{\mathrm{ell}}$ be the uncompactified moduli stack of elliptic curves. Over $\mathcal{M}_{\mathrm{ell}}$, the universal elliptic curve $p\colon \mathcal{E} \to \mathcal{M}_{\mathrm{ell}}$ is defined. This is a stack in its own right, whose fiber over a point $f\colon \mathrm{Spec}\, R \to \mathcal{M}_{\mathrm{ell}}$ is a classical elliptic curve over $R$, namely the curve classified by $f$.

**Definition.** Let $\omega$ be the line bundle of invariant differentials on $\mathcal{M}_{\mathrm{ell}}$. For $n \in 2\mathbf{Z}$ and $m \in \mathbf{Z}$, define a line bundle $L_{n,m}$ on $\mathcal{E}$ by

$$L_{n,m} = p^*\omega^{n/2} \otimes \mathcal{O}_{\mathcal{E}}(me),$$

where $\mathcal{O}_{\mathcal{E}}(me)$ denotes the sheaf of functions on $\mathcal{E}$ that have a pole of order at most $m$ at the identity element $e\colon \mathcal{M}_{\mathrm{ell}} \to \mathcal{E}$ and $\omega = e^*\Omega^1_{\mathcal{E}/\mathcal{M}_{\mathrm{ell}}}$. If $n$ is odd, we set $L_{n,m} = 0$. We denote by

$$\mathrm{JF}_{n,m} = H^0(\mathcal{E}, L_{n,m})$$

the abelian group of global sections.

Note that the direct sum $L_{*,*} = \bigoplus_{n,m \geq 0} L_{n,m}$ is canonically a bundle of bigraded commutative rings.

**Theorem 1.5** (Theorem 2.13). *There is an isomorphism of bigraded rings*

$$\mathrm{JF}_{n,m} \otimes \mathbf{C} \to \mathrm{JF}^{\mathbf{C}}_{n,m}.$$

Having identified Jacobi forms as sections of a line bundle, it is natural to consider the higher cohomology groups as well:

**Definition.** The ring of derived, weakly holomorphic Jacobi forms is a trigraded ring which in tridegree $(n, s, m)$ is given by

$$\mathrm{DJF}_{n,s,m} = H^s(\mathcal{E}, L_{n,m}).$$

This definition is in analogy to derived modular forms:

**Definition.** The ring of derived, weakly holomorphic modular forms is the bigraded ring which in bidegree $(n, s)$ is given by

$$\mathrm{DMF}_{n,s} = H^s(\mathcal{M}_{\mathrm{ell}}, \mathcal{O}_{\mathcal{M}_{\mathrm{ell}}}).$$

Derived Jacobi forms $\mathrm{DJF}_{***}$ form a commutative algebra over $\mathrm{DMF}_{**}$. Unlike in the algebraic case, $\mathrm{DJF}_{n,s,0}$ is not quite isomorphic to $\mathrm{DMF}_{n,s}$ since $H^s(\mathcal{E}; \mathcal{O}_{\mathcal{E}}) \neq H^s(\mathcal{M}_{\mathrm{ell}}; \mathcal{O}_{\mathcal{M}_{\mathrm{ell}}})$ for $s > 0$ in general.

Goerss, Hopkins, and Miller [DFHH14] constructed a sheaf $\mathcal{O}^{\text{top}}$ of $E_\infty$-ring spectra on the stack $\mathcal{M}_{\text{ell}}$, whose global sections are the periodic spectrum of topological modular forms TMF. Moreover, Lurie constructed in [Lur18] a compatible sheaf of $E_\infty$-ring spectra $\mathcal{O}^{\text{top}}_{\mathcal{E}}$ on $\mathcal{E}$.

We can mimick the definition of $L_{0,m} = \mathcal{O}_{\mathcal{E}}(me)$ in spectral algebraic geometry to define sheaves $L_m^{\text{top}}$ on $\mathcal{E}$ with the following properties:

**Theorem 1.6** (Section 3.1)**.**
  (1) *For every $m \geq 0$, there exists a sheaf $L_m^{\text{top}}$ of $\mathcal{O}^{\text{top}}_{\mathcal{E}}$-module spectra on $\mathcal{E}$ such that $\pi_n(L_m^{\text{top}}) \cong L_{n,m}$.*
  (2) *The sum $L_*^{\text{top}} = \bigoplus_{m \geq 0} L_m^{\text{top}}$ has the structure of an $E_2$-algebra spectrum, and the isomorphism $\pi_*(L_*^{\text{top}}) \cong L_{*,*}$ above is an isomorphism of bundles of commutative rings.*
  (3) *There exists a map $a \colon L_m^{\text{top}} \to L_{m+1}^{\text{top}}$ of $\mathcal{O}^{\text{top}}_{\mathcal{E}}$-modules, which on global sections induces the map $L_{m,n} \to L_{m+1,n}$ given by the inclusion $\mathcal{O}_{\mathcal{E}}(me) \to \mathcal{O}_{\mathcal{E}}((m+1)e)$. The colimit*
$$L_\infty^{\text{top}} = \text{colim}(L_0^{\text{top}} \xrightarrow{a} L_1^{\text{top}} \xrightarrow{a} \cdots)$$
*has the structure of a bundle of $E_\infty$-algebra spectra over $\mathcal{O}^{\text{top}}_{\mathcal{E}}$.*

We make the following obvious definition:

**Definition.** The TMF-module spectrum of *topological (weakly holomorphic) Jacobi forms of index $m$* is defined as the global sections
$$\text{TJF}_m = \Gamma_{\mathcal{E}}(L_m^{\text{top}}).$$

It follows from the above theorem that $\text{TJF}_*$ is a graded $E_2$-spectrum, and that $\text{TJF}_\infty$ is an $E_\infty$-ring spectrum. Like TMF, these spectra come with descent spectral sequences
$$E_2^{s,t} = \text{DJF}_{n,s,m} \Rightarrow \pi_n(\text{TJF}_m) \quad \text{with } n = t - s.$$

The computational bulk of this paper is concerned with computing $\text{DJF}_{*,*,*}$ and the coefficients $\pi_*(\text{TJF}_m)$. This is simplified by the following, somewhat surprising, result:

**Theorem 1.7.** *Let $P^m = \text{cofib}(\Sigma \mathbf{C} P^{m-1} \to S^0)$ be the cofiber of the reduced $S^1$-transfer map (define $P^0 = S^0 \vee S^1$ by convention). Then there exists an equivalence of* TMF*-module spectra*
$$\text{TJF}_m \simeq \text{TMF} \wedge P^m.$$
*Under this isomorphism, the map $a \colon L_m^{\text{top}} \to L_{m+1}^{\text{top}}$ corresponds to the map induced by the map $P^m \to P^{m+1}$ coming from the skeletal inclusion $\mathbf{C}P^{m-1} \to \mathbf{C}P^m$. In particular,*
$$\text{TJF}_\infty = \text{colim}_m \text{TJF}_m = (a^{-1} \text{TJF}_*)_0 \simeq \text{TMF} \wedge P^\infty.$$

*Remark* 1.8. This theorem gives an obvious candidate for a connective spectrum $\text{tjF}_m$ of topological weak Jacobi forms, viz., $\text{tmf} \wedge P^m$. However, we have not been able to show that $\text{tjF}_*$ admits the structure of a (graded) ring spectrum, or even that $\text{tjF}_\infty$ is a ring spectrum. Elementary homology considerations easily show that $P^\infty$ on its own does not admit a ring structure.

In the rest of the section, we present our computational results as proven in Sections 4 to 6. The first of these concerns the stable case $\mathrm{TJF}_\infty$. Recall the rings

$$\mathrm{mf}_* = \mathbf{Z}[c_4, c_6, \Delta]/(c_4^3 - c_6^2 - 1728\Delta) \quad \text{and} \quad \mathrm{MF}_* = \mathrm{mf}_*[\Delta^{-1}]$$

of holomorphic and weakly holomorphic integral modular forms. Moreover, the ring $\mathrm{dmf}_{*,*} = H^*(\mathcal{M}_{\mathrm{Weier}}; \omega^{\otimes *})$ of derived modular forms is generated by $c_4$, $c_6$, $\Delta$ in $\mathrm{dmf}_{*,0}$ along with the following classes with the given bidegree (cf. [Bau08]):

$$\begin{array}{lll} |\eta| = (1,1) & |\nu| = (3,1) & |\delta| = (5,1) \\ |\epsilon| = (8,2) & |\beta| = (10,2) & |\kappa| = (14,2) \\ |\bar{\kappa}| = (20,4). & & \end{array}$$

We obtain the ring of derived weakly holomorphic modular forms as

$$\mathrm{DMF}_{**} \cong \mathrm{dmf}_{**}[\Delta^{-1}].$$

**Theorem 1.9.** *We have*

$$\mathrm{DJF}_{*,*,\infty} \cong \mathbf{Z}[b_2, b_3, b_4, b_8, h_1, \Delta^{-1}]/(2h_1, b_3 h_1, 4b_8 + b_4^2 - b_2 b_3^2)$$

*with* $\Delta = -b_2^2 b_8 - 8b_4^3 - 27b_3^4 + 9b_2 b_3^2 b_4$ *in bidegree* $(24, 0)$ *and*

$$|b_i| = (2i, 0); \; h_1 = (1, 1).$$

*where* $|b_i| = 2i$. *The* $(\mathrm{DMF}_{**})$-*module structure is given by*

$$\begin{array}{lll} c_4 \mapsto b_2^2 - 24b_4 & c_6 \mapsto -b_2^3 + 36b_2 b_4 - 216b_3^2 & \Delta \mapsto \Delta \\ \eta \mapsto h_1 & \nu \mapsto 0 & \delta \mapsto b_2 h_1 \\ \epsilon \mapsto 0 & \beta \mapsto 0 & \kappa \mapsto 0 \\ \bar{\kappa} \mapsto 0. & & \end{array}$$

*The descent spectral sequence for* $\mathrm{TJF}_\infty$ *collapses at* $E_4$ *and yields*

$$\pi_*(\mathrm{TJF}_\infty) \cong \mathbf{Z}[\eta, x_2, x_3, x_4, x_4', x_5, x_6, x_8, \Delta^{-1}]/I$$

*with* $|x_i| = 2i$, $|x_4'| = 8$, $|\eta| = 1$, $|\Delta| = 24$,

$$\begin{aligned} I = (&2\eta, \eta^3, x_2\eta, x_3\eta, x_4'\eta, x_5\eta, x_6\eta, \\ &4x_4 - x_2^2, 2x_5 - x_2 x_3, 2x_6 - x_2 x_4', 2x_3 x_4 - x_2 x_5, 2x_4 x_4' - x_2 x_6, \\ &x_5 x_6 - x_3 x_4 x_4', 4x_8 - (x_4')^2 + x_3 x_5). \end{aligned}$$

*and*

$$\Delta = -x_8 x_4 + 9 x_5 x_3 x_4' - 27 x_3^4 - 8(x_4')^3.$$

Note that away from $p = 2$, we have that $\mathrm{DJF}_{*,s,\infty}[\tfrac{1}{2}] = 0$ for $s > 0$. Thus, the descent spectral sequence collapses, and the edge homomorphism

$$\pi_* \mathrm{TJF}_\infty[\tfrac{1}{2}] \to \mathrm{JF}_{*,\infty}[\tfrac{1}{2}]$$

is an isomorphism.

As mentioned before, derived Jacobi forms of index 0 are not just modular forms: Thm. 1.7 says that $\mathrm{TJF}_0 \simeq \mathrm{TMF} \vee \Sigma\, \mathrm{TMF}$, which was already proved in [GM23]. The ring structure on $\mathrm{TJF}_{*,*,0}$ is given as follows:

**Proposition 1.10.** *We have*

$$\pi_* \mathrm{TJF}_0 \cong \pi_* \mathrm{TMF}[\tau]/(\tau^2 - \tau\eta) \quad (|\tau| = 1).$$

We leverage the computation of Thm. 1.9 and the Atiyah-Hirzebruch spectral sequence for $\mathrm{TMF}_*(P^m)$ to obtain:

**Theorem 1.11.** *We have*
$$\mathrm{DJF}_{***}[\tfrac{1}{2}] \cong \mathbf{Z}[\tfrac{1}{2}, a, b_2, b_3, b_4, c_4, c_6, \Delta^{\pm 1}, \tau, \alpha, \beta, \gamma]/I$$
*with the following tridegrees:*
$$|a| = (0,0,1);\ |b_i| = (2i, 0, i);$$
$$|c_4| = (8,0,0);\ |c_6| = (12,0,0);\ |\Delta| = (24,0,0);$$
$$|\tau| = (1,1,0);\ |\alpha| = (3,1,0);\ |\beta| = (10,2,0);\ |\gamma| = (7,1,2),$$
*and where the ideal $I$ is*
$$I = (1728\Delta - c_6^2 + c_4^3, c_4\alpha, c_4\beta, c_4\gamma, c_6\alpha, c_6\beta, c_6\gamma,$$
$$a^4 c_4 - b_2^2 + 24 b_4, a^6 c_6 + b_2^3 - 36 b_2 b_4 + 216 b_3^2,$$
$$a^{12}\Delta - \tfrac{1}{4}b_2^3 b_3^2 + 27 b_3^4 - 9 b_2 b_3^2 b_4 = \tfrac{1}{4}b_2^2 b_4^2 + 8 b_4^3,$$
$$\tau^2, \tau a, \tau b_2, \tau b_3 - a\gamma, \tau b_4, \tau\gamma,$$
$$\alpha^2, 3\alpha, 3\beta, b_2\alpha, b_4\alpha, b_2\beta, b_4\beta, 3\gamma, a^2\gamma, \gamma^2, a^2\beta - \gamma\alpha, b_2\gamma, b_4\gamma, a^2\alpha).$$

*The $(\mathrm{DMF}_{**})_{(3)}$-algebra structure is given by the fact that $\alpha = \nu \in \mathrm{DMF}_{3,1}$ and $\beta \in \mathrm{DMF}_{10,2}$ and the modular forms $c_4$, $c_6$, $\Delta$ map to the classes of the same name in $\mathrm{DJF}$.*

In particular, away from 6, $b_4$ and $\Delta$ can be expressed as polynomials in the remaining generators and all torsion classes become zero, so that we have

**Corollary 1.12.**
$$\mathrm{DJF}_{*,*,*}[\tfrac{1}{6}] \cong \mathrm{JF}_{*,*}[\tfrac{1}{6}, \tau]/(\tau^2, \tau a, \tau b_2, \tau b_3)$$
*with*
$$\mathrm{JF}_{*,*}[\tfrac{1}{6}] = \mathbf{Z}[\tfrac{1}{6}, a, c_4, c_6, b_2, b_3, \Delta^{-1}]/(a^6 c_6 - \tfrac{1}{2}b_2^3 + \tfrac{3}{2}a^4 c_4 b_2 + 216 b_3^2).$$
*In particular, $\mathrm{DJF}_{*,s,m} = 0$ if $s > 0$ and $m > 0$, and the descent spectral sequence collapses, giving*
$$\pi_* \mathrm{TJF}_m[\tfrac{1}{6}] \cong \mathrm{DJF}_{*,0,m} \quad (m > 0).$$

We also compute the descent spectral sequence at $p = 3$. However, the structure of $\pi_*(\mathrm{TJF}_*)_{(3)}$ is not elucidated by writing down generators and relations for this bigraded $\mathrm{TMF}_*$-algebra, but it can be read off of the chart in Figure 5.4.

We leave analogous computations for $p = 2$ for a separate paper. We moreover invite the reader to compare our results on non-derived Jacobi forms with those of [GW20, Theorem 4.1].

## 2. Jacobi forms and Looijenga line bundles

Apart from recalling and establishing some notation and basic facts about Jacobi forms, the goal of this section is twofold: We compute the ring of weak Jacobi forms in Theorem 2.7 and then we use it to identify in Theorem 2.13 Jacobi forms with algebraic sections of Looijenga line bundles (which are even defined over $\mathbf{Z}$). This is an important stepping stone for establishing a viable definition of topological Jacobi forms.

Denote by $\mathrm{mf}_* = \mathbf{Z}[c_4, c_6, \Delta]/(c_4^3 - c_6^2 - 1728\Delta)$ the graded ring of integral modular forms (with the topologically motivated doubling of degrees, $|c_i| = 2i$). Here $c_4$ and $c_6$ can be taken to be suitably scaled Eisenstein series. We also define $\mathrm{MF}_* = \mathrm{mf}[\Delta^{-1}]$ and call its elements weakly holomorphic modular forms; they correspond to meromorphic modular forms (of any integer weight) which are holomorphic away from the cusp.

Recall from Section 1 the definition of the complex groups $\mathrm{jf}^{\mathbf{C}}_{**}$, $\mathrm{jF}^{\mathbf{C}}_{**}$, and $\mathrm{JF}^{\mathbf{C}}_{**}$ of holomorphic, weak, and weakly holomorphic Jacobi forms, respectively. Our weak Jacobi forms agree with those given in [EZ85, p. 104].

The transformation properties (1.1) and (1.2) can be phrased as an invariance property under the action of the holomorph $\Gamma = \mathbf{Z}^2 \ltimes \mathrm{SL}_2(\mathbf{Z})$.

More concretely, given $k \in \mathbf{Z}$ and $m \in \frac{1}{2} \cdot \mathbf{Z}$, define for a function $\phi \colon \mathcal{H} \times \mathbf{C} \to \mathbf{C}$, a matrix $M = \begin{pmatrix} a & b \\ c & d \end{pmatrix} \in \mathrm{SL}_2(\mathbf{Z}) < \Gamma$, and $(\lambda, \mu) \in \mathbf{Z}^2 < \Gamma$:

$$(2.1) \qquad (\phi|_{k,m} M)(\tau, z) = (c\tau + d)^{-k} e\left(\frac{-cmz^2}{c\tau + d}\right) \phi\left(\frac{a\tau + b}{c\tau + d}, \frac{z}{c\tau + d}\right)$$

and

$$(2.2) \qquad (\phi|_{k,m}(\lambda, \mu))(\tau, z) = e\bigl(m(\lambda^2 \tau + 2\lambda z + \lambda + \mu)\bigr) \phi(\tau, z + \lambda\tau + \mu).$$

Then the transformation properties are equivalent with

$$\phi|_{k,m} \gamma = \phi \quad \text{for all } \gamma \in \Gamma.$$

We make a number of immediate observations:

(1) Weakly holomorphic Jacobi forms are 1-periodic in $\tau$. They are also 1-periodic in $z$ up to a sign, which is $+1$ precisely if the index $m$ is an integer. Thus a weakly holomorphic Jacobi form is expressible as a function in $q$ and $\zeta^{\frac{1}{2}}$ even without the Laurent polynomial condition (1.3).
(2) When $m = 0$, condition (1.2) implies that $\phi$ is also $\tau$-periodic in $z$, and thus constant. Thus $\mathrm{jf}^{\mathbf{C}}_{n,0} = \mathrm{jF}^{\mathbf{C}}_{n,0} \cong \mathrm{mf}_n \otimes \mathbf{C}$ and $\mathrm{JF}^{\mathbf{C}}_{n,0} \cong \mathrm{MF}_n \otimes \mathbf{C}$ for all $n \in \mathbf{Z}$.
(3) The groups $\mathrm{jf}^{\mathbf{C}}_{**}$, $\mathrm{jF}^{\mathbf{C}}_{**}$ and $\mathrm{JF}^{\mathbf{C}}_{**}$ all carry the structure of bigraded $\mathbf{C}$-algebras under addition and multiplication of functions. This is straightforward for $\mathrm{jF}^{\mathbf{C}}$ and $\mathrm{JF}^{\mathbf{C}}$; the corresponding fact for holomorphic Jacobi forms $\mathrm{jf}^{\mathbf{C}}$ is [EZ85, Theorem 1.5]. In particular, by the previous point, all three versions of Jacobi forms are algebras over $\mathrm{mf}_*$, and $\mathrm{JF}^{\mathbf{C}}_{**}$ is even an algebra over $\mathrm{MF}_*$.

Using the $\mathrm{mf}_*$-module structure, we have:

**Lemma 2.3.** *The inclusion $\mathrm{jF}^{\mathbf{C}}_{**} \subset \mathrm{JF}^{\mathbf{C}}_{**}$ induces an isomorphism*

$$\mathrm{JF}^{\mathbf{C}}_{**} \cong \mathrm{jF}^{\mathbf{C}}_{**}[\Delta^{-1}].$$

*Proof.* Since $\Delta \in \mathrm{JF}^{\mathbf{C}}_{24,0}$ is invertible in $\mathrm{JF}^{\mathbf{C}}_{**}$, we obtain an inclusion $\mathrm{jF}^{\mathbf{C}}_{**}[\Delta^{-1}] \to \mathrm{JF}^{\mathbf{C}}_{**}$. To show surjectivity, let $\phi \in \mathrm{JF}_{m,n}$ be a weakly holomorphic Jacobi form, and write

$$\phi(\tau, z) = \sum_{n \geq -N} \sum_r c(n, r) q^n \zeta^r.$$

Since $\Delta$ is a cusp form, its $q$-expansion vanishes at 0, i.e. $\Delta(\tau) \in q\mathbf{Z}[\![q]\!]$, and hence $\Delta^N \phi \in \mathrm{jF}^{\mathbf{C}}_{**}$. □

Similarly to modular forms, basic examples of Jacobi forms are given by *Eisenstein series* of weight $k$ and index $m$:

$$(2.4) \qquad E_{k,m}(\tau, z) = \sum_{\gamma \in U \backslash \Gamma} 1|_{k,m}\gamma \in \mathrm{JF}^{\mathbf{C}}_{2k+4m, 2m}$$

where $U \cong \mathbf{Z} \times \mathbf{Z}$ is the subgroup generated by elements of the form $(0, \mu) \ltimes \pm \begin{pmatrix} 1 & b \\ 0 & 1 \end{pmatrix}$ (which satisfy $1|_{k,m}\gamma = 1$).

The invariance of these series holds by definition, and they converge absolutely for $k \geq 4$. The proof that $E_{k,m}$ satisfies the holomorphicity conditions can be found in [EZ85, §I.2]; we will not make use of this.

We recall some classical special functions. The Jacobi theta function is defined by

$$\vartheta(\tau, z) = \sum_{n \in \mathbf{Z}} q^{\frac{1}{2}n^2} \zeta^n$$

and has a variant

$$\vartheta_1(\tau, z) = -\vartheta_{11}(\tau, z) = -iq^{\frac{1}{8}}\zeta^{\frac{1}{2}}\vartheta(\tau, z + \tfrac{1}{2}\tau + \tfrac{1}{2}) = \sum_{n \in \mathbf{Z}+\frac{1}{2}} i^{2n+2} q^{\frac{1}{2}n^2} \zeta^n.$$

Here and elsewhere, when writing $q^a$ or $\zeta^a$ for a rational number $a$, we always mean the branch given by $e(a\tau)$ and $e(az)$, respectively.

The Dedekind eta function is defined by

$$\eta(\tau) = q^{\frac{1}{24}} \prod_{n \geq 1}(1 - q^n)$$

and satisfies $\eta^{24} = \Delta$.[1]

**Lemma 2.5.** *The function*

$$a(\tau, z) = \frac{\vartheta_1(\tau, z)}{\eta(\tau)^3} \in \mathrm{jF}^{\mathbf{C}}_{0,1}$$

*is a weak Jacobi form of weight $-1$ and index $\frac{1}{2}$ with a simple zero at $z = 0$ and no other zeroes in $\mathbf{C}/(\mathbf{Z} + \tau\mathbf{Z})$.*

*Proof.* Since the only zero of $\eta$ within the open unit disk is at $q = 0$, the function $a$ satisfies the holomorphicity condition of a weak Jacobi form. We verify the transformation properties of $a$. Since $\vartheta_1(\tau, z+1) = -\vartheta_1(\tau, z)$ and

$$\vartheta_1(\tau, z+\tau) = -q^{\frac{1}{2}}\zeta^{-1}\vartheta_1(\tau, z),$$

the function $\vartheta_1$, and hence also $a$, transforms according to (1.2) for $m = \frac{1}{2}$.

For the modular transformation, observe that $\vartheta_1(\tau+1, z) = e(\frac{1}{8})\vartheta_1(\tau, z+\frac{1}{2})$ and $\eta(\tau+1) = e(\frac{1}{24})\eta(\tau)$, so that $a(\tau+1, z) = a(\tau)$. For the transformation $(\tau, z) \mapsto (-\frac{1}{\tau}, \frac{z}{\tau})$, we use the Jacobi identity (cf. [Mum07, §I.9])

$$\vartheta_1(-\tfrac{1}{\tau}, \tfrac{z}{\tau}) = e(-\tfrac{3}{8})\tau^{\frac{1}{2}} e(\tfrac{z^2}{2\tau})\vartheta_1(\tau, z)$$

and

$$\eta(-\tfrac{1}{\tau}) = e(-\tfrac{1}{8})\tau^{\frac{1}{2}}\, \eta(\tau)$$

to obtain

$$a(-\tfrac{1}{\tau}, \tfrac{z}{\tau}) = \tau^{-1} e(\tfrac{1}{2}z^2 \tau^{-1})\, a(\tau, z),$$

---

[1] In the complex-analytic literature, the discriminant is often scaled by a factor of $(2\pi)^{12}$.

showing that $a$ is a Jacobi form of weight $-1$ and index $\frac{1}{2}$. Here $\tau^{\frac{1}{2}}$ denotes the root in the upper half plane. Since $a$ is of weight $-1$, the same is true for the modular form $a(-,0)$, which therefore must vanish. The zero at $z=0$ must be simple and $a$ has no other zeros since a Jacobi form of index $m$ has exactly $2m$ zeros, counting multiplicities, by [EZ85, Theorem 1.2] (applied to $a^2$ to ensure integral index). □

*Remark* 2.6. The form $a$ is a square root of the (weight $-2$, index 1) weak Jacobi form $\tilde{\phi}_{-2,1} = \frac{1}{144\Delta}(E_{6,0}E_{4,1} - E_{4,0}E_{6,1})$ classically considered (e.g. [EZ85, Thm. 9.3]). Indeed, [DMZ12, (4.29)] shows that $\tilde{\phi}_{-2,1} = \frac{\vartheta_1^2}{\eta^6}$.

**Theorem 2.7.**
$$\mathrm{jF}^{\mathbf{C}}_{**} \cong \mathbf{C} \otimes \mathrm{mf}_*[a,b,c]/(432c^2 - (b^3 - 3c_4a^4b + 2c_6a^6))$$
with $|a| = (0,1)$, $|b| = (4,2)$, and $|c| = (6,3)$.

*Proof.* The generators $a$, $b$, $c$ are given as follows:
- The element $a$ is the weak Jacobi form of the same name constructed in Lemma 2.5.
- The element $b$ is the weak Jacobi form denoted by $\tilde{\phi}_{0,1} = \frac{\phi_{12,1}}{\Delta}$ in [EZ85, Thm. 9.3], where $\phi_{12,1} = \frac{1}{144}(E_4^2 E_{4,1} - E_6 E_{6,1})$.
- The element $c$ is given by $c(\tau,z) = \frac{\vartheta_1(\tau, 2z)}{\vartheta_1(\tau, z)}$. By [DMZ12, (4.31)], this Jacobi form satisfies $ac = \tilde{\phi}_{-1,2}$ for the class
$$\tilde{\phi}_{-1,2} = \frac{\phi_{11,2}}{\Delta} = \frac{1}{288\pi i \Delta}(E'_{4,1}E_{6,1} - E_{4,1}E'_{6,1})$$
defined in [EZ85, Ch. 9]. (The derivatives are with respect to the $z$ variable.)

In [EZ85, Thm. 9.4], Eichler and Zagier determine the structure of $\mathrm{jF}^{\mathbf{C}}_{*,2*}$, the integer index weak Jacobi forms, as
$$(2.8) \qquad \mathrm{jF}^{\mathbf{C}}_{n,2m} = \mathrm{mf}_*[a^2,b,ac]/(432(ac)^2 - a^2(b^3 - 3c_4a^4b + 2c_6a^6)).$$
Since multiplication by $a$ on $\mathrm{jF}^{\mathbf{C}}_{**}$ is injective, we observe that
$$h\colon \mathrm{mf}_*[a,b,c]/(432c^2 - (b^3 - 3c_4a^4b + 2c_6a^6)) \to \mathrm{jF}_{*,*}$$
is indeed injective. Now let $f \in \mathrm{jF}_{*,*}$ be of half-integer index. Then $af \in \mathrm{jF}_{*,*}$ is a polynomial in $a^2, b$ and $ac$ by (2.8). To see that $f$ is in the image of $h$, we must show that a polynomial in $a^2, b$ and $ac$ is not divisible by $a$ (as a Jacobi form) if there is a monomial of the form $b^n$. But indeed, $b$ has a Taylor expansion of the form $12 + O(z)$ by Formula (19) in Chapter 3 of [EZ85]). Since $a$ is zero at $z=0$, we see that $b^n$ is indeed not divisible by $a$. Thus, $h$ is surjective as well. □

Multiplication by $a$ gives an injection
$$a\colon \mathrm{jF}^{\mathbf{C}}_{n,m} \to \mathrm{jF}^{\mathbf{C}}_{n,m+1}.$$
Denote by $\mathrm{jF}^{\mathbf{C}}_{\infty,n}$ the (singly-graded) colimit
$$\mathrm{jF}^{\mathbf{C}}_{n,\infty} = \mathrm{colim}(\mathrm{jF}^{\mathbf{C}}_{n,0} \xrightarrow{a} \mathrm{jF}^{\mathbf{C}}_{n,1} \xrightarrow{a} \cdots)$$
We immediately get that
$$\mathrm{jF}^{\mathbf{C}}_{*,\infty} \cong \mathrm{mf}_*[b,c]/(432c^2 - b^3 - 3c_4b + 2c_6)$$
with $|b|=4$, $|c|=6$.

In particular, $jF^{\mathbf{C}}_{n,m}$ is finite-dimensional for each $(n,m)$ (and hence $jf^{\mathbf{C}}$ as well, cf. [EZ85, Theorem 1.1]). However, while many partial results on $jf^{\mathbf{C}}_{n,m}$ are known, no explicit formula for the ring structure such as Thm. 2.7 is. In general, even the dimension of $jf^{\mathbf{C}}_{n,m}$ is hard to compute (cf. [EZ85, Section 10], [Sko07]).

**The Looijenga line bundle.** The Looijenga line bundle $L^{\mathbf{C}}_{m,n}$ is a complex line bundle on $\mathcal{E} \otimes \mathbf{C}$, the universal complex elliptic curve living over the complexified moduli stack of elliptic curves $\mathcal{M}_{\text{ell}} \otimes \mathbf{C}$. We will show that the sections of $L^{\mathbf{C}}_{m,n}$ are precisely weakly holomorphic Jacobi forms of index $m$ and dimension $n$.

While the general theory of Looijenga line bundles [Loo77, Rez20] associated to quadratic forms on lattices is much richer, we will only look at the one-dimensional case, in which every quadratic form $Q\colon \mathbf{Z} \to \mathbf{Z}$ is uniquely determined by $Q(1) = 2m$, where $m$ corresponds to the index of the Jacobi form. Here, unlike in Looijenga's definition, $m$ is allowed to be a half integer.

To define the complex Looijenga line bundle, consider the following commutative diagram of stacks:

$$\begin{array}{ccc} \mathcal{H} \times \mathbf{C} & \xrightarrow{\pi_{\mathcal{E}}} & \mathcal{E} \otimes \mathbf{C} \\ \downarrow{p_1} & & \downarrow{p} \\ \mathcal{H} & \xrightarrow{\pi_{\mathcal{M}_{\text{ell}}}} & \mathcal{M}_{\text{ell}} \otimes \mathbf{C} \end{array}$$

where the horizontal maps $\pi_{\mathcal{E}}$ and $\pi_{\mathcal{M}_{\text{ell}}}$ are stack covers with respect to the group actions of $\Gamma = \mathbf{Z}^2 \rtimes \mathrm{SL}_2(\mathbf{Z})$ and $\mathrm{SL}_2(\mathbf{Z})$, respectively, given by

$$\left(\begin{smallmatrix} a & b \\ c & d \end{smallmatrix}\right).(\tau, z) = \left(\frac{a\tau + b}{c\tau + d}, \frac{z}{c\tau + d}\right) \quad \text{for } \left(\begin{smallmatrix} a & b \\ c & d \end{smallmatrix}\right) \in \mathrm{SL}_2(\mathbf{Z}) < \Gamma$$

$$(\lambda, \mu).(\tau, z) = (\tau, z + \lambda\tau + \mu) \quad \text{for } (\lambda, \mu) \in \mathbf{Z}^2 < \Gamma.$$

**Definition.** For $m \in \frac{1}{2}\mathbf{Z}$ and $n \in 2\mathbf{Z}$, the *complex Looijenga line bundle* $L^{\mathbf{C}}_{n,2m}$ on $\mathcal{E} \otimes \mathbf{C}$ is uniquely defined (by descent) by the properties

(1) $\pi_{\mathcal{E}}^*(L^{\mathbf{C}}_{n,2m}) \cong (\mathcal{H} \times \mathbf{C}) \times \mathbf{C}$ is a trivial line bundle, and
(2) The monodromy action of $\Gamma$ on $\pi_{\mathcal{E}}^*(L^{\mathbf{C}}_{n,2m})$ is given by

$$\gamma.(v, \tau, z) = \Big((\phi|_{k,m}\gamma)v, \gamma.(\tau, z)\Big)$$

with $\gamma \in \Gamma$, $k = \frac{n-4m}{2}$ and $\phi|_{k,m}$ as in Equations (2.1) and (2.2).

**Lemma 2.9.** *There is an isomorphism $L^{\mathbf{C}}_{n+2,m} \cong L^{\mathbf{C}}_{n,m} \otimes p^*(\omega)$.*

*Proof.* Recall that the line bundle $\pi_{\mathcal{M}_{\text{ell}}}^*(\omega^{\mathbf{C}}) \cong \mathbf{C} \times \mathcal{H}$ is trivial with $\mathrm{SL}_2(\mathbf{Z})$-action given by

$$\left(\begin{smallmatrix} a & b \\ c & d \end{smallmatrix}\right).(v, \tau) = \Big((c\tau + d)^{-1}v, \frac{a\tau + b}{c\tau + d}\Big).$$

Since

$$\frac{\phi|_{k+1,m}\gamma}{\phi|_{k,m}\gamma} = \begin{cases} (c\tau+d)^{-1}; & \gamma = \left(\begin{smallmatrix} a & b \\ c & d \end{smallmatrix}\right) \\ 1; & \gamma = (\lambda, \mu) \in \mathbf{Z}^2, \end{cases}$$

the claim follows. $\square$

By construction, $L^{\mathbf{C}}_{m,n}$ is a holomorphic line bundle on $\mathcal{E} \otimes \mathbf{C}$. We now construct an isomorphism to an algebraic line bundle, which is even defined over $\mathbf{Z}$.

**Proposition 2.10.** *There is an isomorphism $L^{\mathbf{C}}_{2n,m} \cong p^*\omega^n \otimes \mathcal{O}_{\mathcal{E}}(me) \otimes \mathbf{C}$.*

*Proof.* By Lemma 2.5, $a$ has a simple zero at $z = 0$ and no other zeros in $\mathbf{C}/(\mathbf{Z}+\tau\mathbf{Z})$. Therefore, it induces a trivialization of

$$L_{0,1}^{\mathbf{C}} \otimes \mathcal{O}(-e)$$

and thus, using Lemma 2.9, an isomorphism of line bundles

$$L_{2n,m}^{\mathbf{C}} \cong \mathcal{O}(me) \otimes p^*\omega^n \otimes \mathbf{C}. \qquad \square$$

By construction, the diagram

$$\begin{array}{ccc}
\mathcal{O}_{\mathcal{E}}(me) \otimes \omega^n & \longrightarrow & \mathcal{O}_{\mathcal{E}}((m+1)e) \otimes \omega^n \\
\downarrow a^m & & \downarrow a^{m+1} \\
L_{2n,m}^{\mathbf{C}} & \xrightarrow{a} & L_{2n,m+1}^{\mathbf{C}}
\end{array}$$

commutes. Note that the vertical maps are isomorphisms by the proof of the preceding proposition.

**An algebraic theory of Jacobi forms.** In this section, we will prove an analogue of the classic theorem that the space of weakly holomorphic modular forms of weight $k$ is isomorphic to that of sections of $\omega^k$ on $\mathcal{M}_{\mathrm{ell}} \otimes \mathbf{C}$ (see [DR73, Section VII.4] and [MO20, Section A.2.3]).

**Lemma 2.11.** *Let $s$ be an (algebraic) section of $\mathcal{O}_{\mathcal{E}}(me) \otimes \omega^k \otimes \mathbf{C}$ over $\mathcal{M}_{\mathrm{ell}} \otimes \mathbf{C}$. Under the isomorphism from Proposition 2.10, this corresponds to a section of $L_{m,2k}^{\mathbf{C}}$ and thus to a holomorphic function $f$ on $\mathcal{H} \times \mathbf{C}$ satisfying the transformation law of Jacobi forms. This $f$ is a weak Jacobi form, i.e. its Fourier expansion is of the form $\sum_{n \geq -M, r \in \mathbf{Z}} c(n,r) q^n \zeta^r$.*

*Proof.* As for the analogous result for modular forms, we utilize the Tate curve. Let $\mathsf{Conv} \subset \mathbf{C}((q))$ be the subring of those Laurent series defining a holomorphic function on the punctured open unit disk $\mathbb{D}^{\circ}$ and let $\mathrm{ev}_q \colon \mathsf{Conv} \to \mathbf{C}$ be the evaluation at $q \in \mathbb{D}^{\circ}$. The Tate curve $\mathsf{Tate}$ is an elliptic curve over $\mathsf{Conv}$ such that $(\mathrm{ev}_q)_* \mathsf{Tate} \cong \mathbf{C}^{\times}/q^{\mathbf{Z}}$. It can be written in Weierstrass form [Sil94, Theorem V.1.1]

$$(2.12) \qquad y^2 + xy = x^3 + a_4(q)x + a_6(q) \quad (a_4, a_6 \in \mathsf{Conv})$$

and thus its bundle of invariant differentials on $\mathsf{Conv}$ is trivial (see [MO20, Appendix A.3]).

We can pull back $s$ along the classifying map $\mathsf{Conv} \to \mathcal{M}_{\mathrm{ell}} \otimes \mathbf{C}$ of $\mathsf{Tate}$ to obtain a section $\bar{s}$ of $\mathcal{O}_{\mathsf{Tate}}(me) \otimes \omega^k \cong \mathcal{O}_{\mathsf{Tate}}(me)$. As on every Weierstrass elliptic curve, the coordinate functions $x$ and $y$ from (2.12) generate $\bigoplus_m \Gamma(\mathcal{O}_{\mathsf{Tate}}(me))$ as a $\mathsf{Conv}$-algebra. Since for every $q$, the pushforward $(\mathrm{ev}_q)_* \mathsf{Tate}$ is a quotient of $\mathbf{C}^{\times}$, we can view $x$ and $y$ as meromorphic functions in $q \in \mathbb{D}$ and $\zeta \in \mathbf{C}^{\times}$. As in the proof of [Sil94, Theorem V.3.1(c)] (second and third formula on p.426 loc.cit.), these can be written in the form $\sum_{n \geq 0, r} d(n,r) q^n \zeta^r$. Thus, $\bar{s}$ can be written in the form $\sum_{n \geq -M, r} c(n,r) q^n \zeta^r$.

We obtain the corresponding section of $L_{n,m}^{\mathbf{C}}$ by multiplying by $a^{\otimes m}$. Since $a$ is a weak Jacobi form, the result can be written in the same form (with different $c(n,r)$) when viewed as a function on $\mathbb{D}^{\circ} \times \mathbf{C}^{\times}$. (Here, we use implicitly that $L_{n,m}^{\mathbf{C}}$ is by construction trivial on $\mathbb{D}^{\circ} \times \mathbf{C}^{\times}$.) The result follows. $\qquad \square$

**Theorem 2.13** (Theorem 1.5). *For all $m, k$, we have*
$$\mathrm{JF}^{\mathbf{C}}_{2k,m} \cong \Gamma_{\mathcal{E}\otimes\mathbf{C}}(\mathcal{O}_{\mathcal{E}}(me) \otimes \omega^k \otimes \mathbf{C}).$$
*These isomorphisms define an isomorphism of bigraded rings.*

*Proof.* By definition, holomorphic sections of $L^{\mathbf{C}}_{2k,m}$ correspond to holomorphic $\Gamma$-invariant functions on $\mathcal{E} \otimes \mathbf{C}$, i.e. functions satisfying the transformation law for Jacobi forms. It remains to show that a section of $L^{\mathbf{C}}_{2k,m} \cong \mathcal{O}_{\mathcal{E}}(me) \otimes \omega^k \otimes \mathbf{C}$ is algebraic if and only if the Fourier expansion of the corresponding function is a Laurent series of the form
$$\sum_{n \geq -M, r} c(n,r) q^n \zeta^r;$$
we will call such sections *Laurent*. For ease of notation, we will omit the $\otimes \mathbf{C}$ in the rest of the proof.

By the preceding lemma, every algebraic section is Laurent; thus, multiplication with $a^m$ defines a natural homomorphism $\Gamma_{\mathcal{E}}(\mathcal{O}_{\mathcal{E}}(me) \otimes \omega^k) \to \mathrm{JF}_{2k,m}$. By Thm. 2.7, we have a short exact sequence
$$0 \to \mathrm{JF}_{2k,m} \xrightarrow{\cdot a} \mathrm{JF}_{2k,m+1} \to \begin{cases} 0; & m=0 \\ \mathrm{MF}_{2k-2m-2}; & m > 0 \end{cases} \to 0$$

For $m > 0$, a section $s \colon \mathrm{MF}_{2k-2m-2} \to \mathrm{JF}_{2k,m+1}$ of this sequence is given by
$$s(\phi) = \begin{cases} b^{\frac{m+1}{2}} \phi; & m \text{ odd} \\ cb^{\frac{m-2}{2}} \phi; & m \geq 2 \text{ even} \end{cases}$$

By [Del75, Section 1], we also have a short exact sequence
$$0 \to p_*\mathcal{O}_{\mathcal{E}}(me) \otimes \omega^k \to p_*\mathcal{O}_{\mathcal{E}}((m+1)e) \otimes \omega^k \to \begin{cases} 0; & m = 0 \\ \omega^{k-m-1}; & m > 0 \end{cases} \to 0$$

and the (algebraic) sections of the quotient are $\mathrm{MF}_{2k-2m-2}$. Concretely, we obtain this modular form from a section $p_*\mathcal{O}_{\mathcal{E}}((m+1)e) \otimes \omega^n$ by Taylor expanding and taking the lowest coefficient (of $z^{-m-1}$).

By [EZ85, Formula (19), p.40], $b = \frac{\phi_{12,1}}{\Delta} = 12 + O(z^2)$. Similarly, $c = \frac{\vartheta_1(q,2z)}{\vartheta_1(q,z)} = 2 + O(z)$, so the lowest term in the $z$-Taylor expansion of $c^i b^j$ is constant in $\tau$. Using [EZ85, Formula (19), p.40] again and the identity $a^2 = \frac{\phi_{10,1}}{\Delta}$ from Remark 2.6, we also obtain $a = \alpha z + O(z^2)$ for a nonzero constant $\alpha$. Thus, the map $\mathrm{JF}_{2k,m+1} \to \mathrm{MF}_{2k-2m-2}$ is, up to a constant only dependent on the index, the constant term of $z$-Taylor expansion. Thus we obtain in total a commutative diagram (for $m > 0$)

$$\begin{array}{ccccccccc}
0 & \longrightarrow & \Gamma_{\mathcal{E}}(\mathcal{O}_{\mathcal{E}}(me) \otimes \omega^k) & \longrightarrow & \Gamma_{\mathcal{E}}(\mathcal{O}_{\mathcal{E}}((m+1)e) \otimes \omega^k) & \longrightarrow & \mathrm{MF}_{2k-2m-2} & \longrightarrow & 0 \\
& & \downarrow a^m & & \downarrow a^{m+1} & & \downarrow \cong & & \\
0 & \longrightarrow & \mathrm{JF}_{2k,m} & \xrightarrow{a} & \mathrm{JF}_{2k,m+1} & \longrightarrow & \mathrm{MF}_{2k-2m-2}, & \longrightarrow & 0
\end{array}$$

where the rightmost vertical map is multiplication by a suitable constant. This shows that the vertical arrow $a^{m+1}$ is an isomorphism if $a^m$ is an isomorphism (for

$m > 0$). To ground the induction, observe that in the diagram

$$\begin{array}{ccc}
\Gamma_{\mathcal{E}}(\mathcal{O}_{\mathcal{E}} \otimes \omega^k) & \longrightarrow & \Gamma_{\mathcal{E}}(\mathcal{O}_{\mathcal{E}}(e) \otimes \omega^k) \\
\downarrow & & \downarrow \\
\mathrm{JF}_{2k,0} & \xrightarrow{a} & \mathrm{JF}_{2k,1}
\end{array}$$

the horizontal arrows are isomorphisms (using Riemann–Roch in the top row and Theorem 2.7 in the bottom row) and the left vertical arrow is also an isomorphism (as both source and target are just modular forms of weight $2k$). Thus, the right vertical arrow is an isomorphism and induction yields our result. □

This theorem allows us to set up an integral version of Jacobi forms.

**Definition.** We set
$$L_{2n,m} = p^*\omega^n \otimes \mathcal{O}_{\mathcal{E}}(me)$$
and
$$\mathrm{JF}_{2n,m} = \Gamma_{\mathcal{E}}(L_{2n,m}).$$

By the preceding theorem, we have indeed $L_{n,m} \otimes \mathbf{C} \cong L_{n,m}^{\mathbf{C}}$ and $\mathrm{JF}_{n,m} \otimes \mathbf{C} \cong \mathrm{JF}_{n,m}^{\mathbf{C}}$.

*Remark* 2.14. In [Kra95], Kramer sets up an integral version of Jacobi forms as well, actually also for higher-dimensional variants. Restricted to our one-dimensional case, the theory is based on the pullback of the Poincaré line bundle on $E \times E^\vee$ along the graph $E \to E \times E^\vee$ of an autoduality $E \xrightarrow{\cong} E^\vee$. One can show that this pullback is (with suitable choices) indeed isomorphic to $\mathcal{O}_{\mathcal{E}}(2e) \otimes \omega^2$, giving the equivalence of the two approaches in the integer index case.

3. IDENTIFYING $\mathrm{TJF}_m$ VIA EQUIVARIANT TOPOLOGICAL MODULAR FORMS

In this section, we identify topological Jacobi forms $\mathrm{TJF}_m$ in terms of stunted projective spaces by interpreting $\mathrm{TJF}_m$ in terms of $T$-equivariant TMF, where $T \cong S^1$ is the circle group. We start by recalling its basic properties.

3.1. **Equivariant topological modular forms.** For an even-periodic complex oriented cohomology theory $E$, the ring $E^0(BT)$ defines a formal group law over $E_0$. Viewed differently, $E^0(BT)$ is the Borel $T$-equivariant cohomology of a point. One may wish to *decomplete* this by replacing Borel cohomology by a *genuine* cohomology theory. For instance, the formal group associated to complex $K$-theory is the multiplicative formal group, while $KU_T^0(\{*\}) = \mathrm{Rep}(T) \cong \mathbf{Z}[t^{\pm 1}]$ corepresents the actual multiplicative group. The same idea has been used by Grojnowski [Gro07] for elliptic cohomology theories: while the formal group law is the completion of the elliptic curve, $T$-equivariant elliptic cohomology allows one to recover the elliptic curve itself, essentially as the $T$-equivariant cohomology of a point.

Topological modular forms are, of course, not complex oriented, but they are built out of (complex oriented) elliptic cohomology theories. Classically, TMF was constructed by building a sheaf of $E_\infty$-ring spectra $\mathcal{O}^{\mathrm{top}}$ on $\mathcal{M}_{\mathrm{ell}}$ with the property that $\pi_{2k}\mathcal{O}^{\mathrm{top}} \cong \omega^k$ and whose fibers were elliptic cohomology theories. Periodic TMF was then defined as the global sections of this sheaf, an $E_\infty$-ring spectrum in its own right.

To obtain a *T*-equivariant version of TMF, this definition needed to be reinterpreted. In [Lur09, Lur18], Lurie gave an alternative construction of TMF entirely in terms of spectral algebraic geometry. We no longer define a sheaf of $E_\infty$-ring spectra on an object with an algebraic universal property, but the category of stacks is replaced with the $\infty$-category of (nonconnective) spectral Deligne-Mumford stacks. Lurie constructs a spectral stack $\mathcal{M}_{\text{ell}}^{\text{or}}$ of *oriented* spectral elliptic curves, together with a universal oriented curve $\mathcal{E}^{\text{or}} \to \mathcal{M}_{\text{ell}}^{\text{or}}$. The underlying ordinary stack of $\mathcal{M}_{\text{ell}}^{\text{or}}$ is the classical Deligne-Mumford stack $\mathcal{M}_{\text{ell}}$, and $\mathcal{O}^{\text{top}}$ can be identified with the structure sheaf $\mathcal{O}_{\mathcal{M}_{\text{ell}}^{\text{or}}}$. But $\mathcal{E}^{\text{or}}$ is a spectral Deligne-Mumford stack with its own structure sheaf $\mathcal{O}_{\mathcal{E}^{\text{or}}}$ as well, which can be accessed using the above ideas of *T*-equivariant cohomology. Based on this work of Lurie and older insights of Grojnowski, Gepner and the second author construct in [GM23] a *T*-equivariant version of TMF-cohomology with the following properties:

**Theorem 3.1** (Gepner-Meier)**.** *There exists a colimit-preserving functor of $\infty$-categories*

$$\mathcal{T}\text{MF}_T \colon \{\text{pointed } T\text{-complexes}\} \to \{\text{quasicoherent } \mathcal{O}_{\mathcal{E}^{\text{or}}}\text{-module sheaves}\}^{\text{op}}$$

*with the following properties:*

(1) $\mathcal{T}\text{MF}_T(S^0) \simeq \mathcal{O}_{\mathcal{E}^{\text{or}}}$;
(2) $\mathcal{T}\text{MF}_T(T_+) \simeq e_*\mathcal{O}_{\mathcal{M}_{\text{ell}}^{\text{or}}}$, *for* $e \colon \mathcal{M}_{\text{ell}}^{\text{or}} \to \mathcal{E}^{\text{or}}$ *the unit section;*
(3) $\mathcal{T}\text{MF}_T$ *is symmetric monoidal on finite complexes; in particular, for finite pointed T-complexes X and Y, there is a natural "Künneth" isomorphism*

$$\mathcal{T}\text{MF}_T(X) \otimes_{\mathcal{O}_{\mathcal{E}^{\text{or}}}} \mathcal{T}\text{MF}_T(Y) \xrightarrow{\simeq} \mathcal{T}\text{MF}_T(X \wedge Y);$$

(4) $\mathcal{T}\text{MF}_T$ *factors via $\Sigma^\infty$ through a functor*

$$\{\text{genuine } T\text{-spectra}\} \to \{\text{quasicoherent } \mathcal{O}_{\mathcal{E}^{\text{or}}}\text{-module sheaves}\}^{\text{op}},$$

*which is symmetric monoidal on finite T-spectra and which we also denote by $\mathcal{T}\text{MF}_T$;*
(5) $\text{TMF}_T =_{\text{def}} \Gamma \mathcal{T}\text{MF}_T \colon \{\text{pointed } T\text{-complexes}\} \to \{\text{Spectra}\}^{\text{op}}$ *is a T-equivariant cohomology theory and is thus represented by a genuine T-spectrum, which we denote by the same symbol* $\text{TMF}_T$. *Its underlying nonequivariant spectrum is* TMF.

*Proof.* For every oriented spectral elliptic curve $E$, [GM23] constructs a functor $\widetilde{\mathcal{E}ll}_T$ from finite pointed *T*-complex to {quasicoherent $\mathcal{O}_{\mathcal{E}^{\text{or}}}$-module sheaves}$^{\text{op}}$. The functor $\mathcal{T}\text{MF}_T$ is obtained by specializing $E$ to be the universal oriented spectral elliptic curve $\mathcal{E}^{\text{or}}$ and extending $\mathcal{T}\text{MF}_T = \widetilde{Ell}_T$ by filtered colimits to be defined on all pointed *T*-complexes. The first two properties hold by the construction, and the other three are Theorem 8.1, Proposition 9.2 and Construction 9.3 in [GM23]. □

For a *T*-representation $V$, denote by $S^V$ its one-point compactification. We denote by $S^n$ the sphere with the trivial action and by $T$ the one-sphere with its tautological *T*-action. Consider the cofiber sequence of pointed *T*-complexes

$$T_+ \to S^0 \to S^\rho,$$

where $\rho$ is the tautological 1-dimensional complex representation of $T$. Since $\mathcal{T}\text{MF}_T$ is colimit-preserving, this gives rise to a cofiber sequence

(3.2) $$\mathcal{T}\text{MF}_T(S^\rho) \to \mathcal{O}_{\mathcal{E}^{\text{or}}} \to e_*(\mathcal{O}_{\mathcal{M}_{\text{ell}}^{\text{or}}}).$$

We define $\mathcal{O}_{\mathcal{E}^{\mathrm{or}}}(-e) = \mathcal{T}\mathrm{MF}_\mathrm{T}(S^\rho)$ and, more generally,
$$\mathcal{O}_{\mathcal{E}^{\mathrm{or}}}(-me) = \mathcal{T}\mathrm{MF}_\mathrm{T}(S^\rho)^{\otimes_{\mathcal{O}_{\mathcal{E}^{\mathrm{or}}}} m} \cong \mathcal{T}\mathrm{MF}_\mathrm{T}((S^\rho)^{\wedge m}) \cong \mathcal{T}\mathrm{MF}_\mathrm{T}(S^{m\rho})$$
for $m \geq 0$. As a special case of [GM23, Lemma 8.1], $\mathcal{O}_{\mathcal{E}^{\mathrm{or}}}(-e)$ is invertible. Thus, we can define $\mathcal{O}_{\mathcal{E}^{\mathrm{or}}}(me) = \mathrm{Hom}_{\mathcal{O}_{\mathcal{E}^{\mathrm{or}}}}(\mathcal{O}_{\mathcal{E}^{\mathrm{or}}}(-me), \mathcal{O}_{\mathcal{E}^{\mathrm{or}}})$.

**Theorem 3.3.** *The graded quasicoherent $\mathcal{O}_{\mathcal{E}^{\mathrm{or}}}$-module*
$$\mathcal{O}_{\mathcal{E}^{\mathrm{or}}}(*) = \bigoplus_{m \in \mathbf{Z}} \mathcal{O}_{\mathcal{E}^{\mathrm{or}}}(me)$$
*has an $E_2$-algebra structure, i.e. there is a an $E_2$-monoidal functor $\mathbf{Z} \to \mathrm{QCoh}(\mathcal{E}^{\mathrm{or}})$ sending $m$ to $\mathcal{O}_{\mathcal{E}^{\mathrm{or}}}(me)$.*

*We have furthermore*
$$\pi_{2k}\mathcal{O}_{\mathcal{E}^{\mathrm{or}}}(me) \cong p^*\omega^k \otimes \mathcal{O}_{\mathcal{E}}(me),$$
*where $p\colon \mathcal{E} \to \mathcal{M}_{\mathrm{ell}}$ denotes the structure map.*

*Proof.* For the first part, we follow an argument in [Lur15]: Let $\mathcal{U}$ be a complete $T$-universe, and let $\bigoplus_\mathbb{N} \rho \cong \mathcal{U}_\rho \subset \mathcal{U}$ be the $\rho$-typical subspace. The space $BU(\mathcal{U})$ has a model as the complex Grassmannian of $\mathcal{U}$. Thus, we have an inclusion
$$\mathbf{C}P(\mathcal{U}_\rho) \subset BU(\mathcal{U})$$
of 1-dimensional subrepresentations isomorphic to $\rho$. The $T$-action on $\mathbf{C}P(\mathcal{U}_\rho)$ is trivial and thus $\Omega^2 \mathbf{C}P(\mathcal{U}_\rho)$ is the discrete $T$-space $\mathbf{Z}$ with trivial action. On the other hand, $T$-equivariant Bott periodicity identifies $\Omega^2 BU(\mathcal{U})$ with $R(T) \times BU(\mathcal{U})$. Postcomposing with the $T$-equivariant $J$-homomorphism, yields an $E_2$-monoidal functor
$$\mathbf{Z} \to R(T) \times BU(\mathcal{U}) \to \mathrm{Sp}_T$$
to genuine $T$-spectra, sending 1 to $\Sigma^\infty S^\rho$ and hence $m$ to $\Sigma^\infty S^{m\rho}$. Dualizing and then applying $\mathcal{T}\mathrm{MF}_\mathrm{T}$ thus sends $m$ to $\mathcal{O}_\mathcal{E}(me)$, as demanded.

Applying homotopy groups to the cofiber sequence (3.2) yields short exact sequences of $\mathcal{O}_\mathcal{E}$-module sheaves
$$0 \to \pi_{2k}\mathcal{O}_\mathcal{E}(-e) \to \pi_{2k}\mathcal{O}_\mathcal{E} \to \pi_{2k}e_*(\mathcal{O}_{\mathcal{M}_{\mathrm{ell}}^{\mathrm{or}}}) \to 0.$$
We have that
$$\pi_{2k}\mathcal{O}_{\mathcal{E}^{\mathrm{or}}} \cong \mathcal{O}_\mathcal{E} \otimes p^*\omega^k \text{ and } \pi_{2k}e_*(\mathcal{O}_{\mathcal{M}_{\mathrm{ell}}^{\mathrm{or}}}) \cong e_*\omega^k \cong (e_*\mathcal{O}_{\mathcal{M}_{\mathrm{ell}}}) \otimes_{\mathcal{O}_\mathcal{E}} p^*\omega^k.$$
Under these identifications, $\pi_{2k}\mathcal{O}_\mathcal{E} \to \pi_{2k}e_*(\mathcal{O}_{\mathcal{M}_{\mathrm{ell}}^{\mathrm{or}}})$ is induced by the canonical map $e^*\colon \mathcal{O}_\mathcal{E} \to e_*\mathcal{O}_{\mathcal{M}_{\mathrm{ell}}}$. Since $\ker(e^*) = \mathcal{O}_\mathcal{E}(-e)$, the result follows for $m = -1$. For other $m$, we first use the symmetric monoidality of $\mathcal{T}\mathrm{MF}_\mathrm{T}$ from Part (3) of Theorem 3.1 to get that $\mathcal{O}_{\mathcal{E}^{\mathrm{or}}}(me) \simeq \mathcal{O}_{\mathcal{E}^{\mathrm{or}}}(-e)^{\otimes -m}$. We further use that on quasicoherent $\mathcal{O}_\mathcal{E}$-modules with flat homotopy groups, the functor
$$\pi_*\colon \mathrm{QCoh}(\mathcal{E}) \to \text{quasi-coherent } \pi_*\mathcal{O}_\mathcal{E}\text{-modules}$$
is symmetric monoidal (by a degenerating Künneth spectral sequence). This implies in particular to invertible sheaves. Thus,
$$\pi_*\mathcal{O}_{\mathcal{E}^{\mathrm{or}}}(me) \cong \pi_*\mathcal{O}_{\mathcal{E}^{\mathrm{or}}}(-e)^{\otimes_{\pi_*\mathcal{O}_{\mathcal{E}^{\mathrm{or}}}}(-m)} \qquad \square.$$

**Proposition 3.4.** *For a finite genuine $T$-spectra $X$, there is a natural equivalence*
$$(\mathrm{TMF}_T \wedge DX)^T \simeq \Sigma D_{\mathrm{TMF}}(\mathrm{TMF}_T \wedge X)^T,$$
*where $D$ denotes the Spanier–Whitehead dual.*

16                TILMAN BAUER AND LENNART MEIER*Proof.* Let $p\colon \mathcal{E}^{\mathrm{or}} \to \mathcal{M}_{\mathrm{ell}}^{\mathrm{or}}$ be the projection. For $\mathcal{F} \in \mathrm{QCoh}(\mathcal{E}^{\mathrm{or}})$, [GKMP, Theorem 6.5] gives a natural equivalence

$$(3.5) \qquad \mathrm{Hom}_{\mathcal{O}_{\mathcal{M}_{\mathrm{ell}}^{\mathrm{or}}}}(p_*\mathcal{F}, \mathcal{O}_{\mathcal{M}_{\mathrm{ell}}^{\mathrm{or}}}) \simeq \Sigma^{-1}\mathrm{Hom}_{\mathcal{O}_{\mathcal{E}^{\mathrm{or}}}}(\mathcal{F}, \mathcal{O}_{\mathcal{E}^{\mathrm{or}}}).$$

Since $\Gamma\colon \mathrm{QCoh}(\mathcal{M}_{\mathrm{ell}}^{\mathrm{or}}) \to \mathrm{Mod}_{\mathrm{TMF}}$ is a symmetric monoidal equivalence by [MM15], the left-hand term can be computed as $\mathrm{Hom}(\Gamma(p_*\mathcal{F}), \mathrm{TMF}) = D_{\mathrm{TMF}}\Gamma(\mathcal{F})$.

We now specialize to $\mathcal{F} = \mathcal{TMF}_T(X)$ for a finite genuine $T$-spectrum $X$. Using that $\mathcal{TMF}_T$ is symmetric monoidal and thus preserves duals, we compute

$$\begin{aligned}
(\mathrm{TMF}_T \wedge DX)^T &\simeq \Gamma\, \mathcal{TMF}_T(X) \\
&\simeq \Gamma\mathcal{H}om_{\mathcal{O}_{\mathcal{E}^{\mathrm{or}}}}(\mathcal{TMF}_T(DX), \mathcal{O}_{\mathcal{E}^{\mathrm{or}}}) \\
&= \mathrm{Hom}_{\mathcal{O}_{\mathcal{E}^{\mathrm{or}}}}(\mathcal{TMF}_T(DX), \mathcal{O}_{\mathcal{E}^{\mathrm{or}}}) \\
&\stackrel{(3.5)}{\simeq} \Sigma D_{\mathrm{TMF}}\Gamma(\mathcal{TMF}_T(DX)) \\
&\simeq \Sigma D_{\mathrm{TMF}}(\mathrm{TMF}_T \wedge X)^T. \qquad \square
\end{aligned}$$

3.2. **The identification of** $\mathrm{TJF}_m$**.** Let us recall and refine the definition of $\mathrm{TJF}_m$ from Section 1.

**Definition.** We define the spectrum $\mathrm{TJF}_m$ of *topological Jacobi forms* as

$$\mathrm{TMF}_T(S^{-m\rho}) = \Gamma\mathcal{O}_{\mathcal{E}}(me).$$

We define $\mathrm{TJF}_\infty$ as $\mathrm{colim}_m \mathrm{TJF}_m$.

*Remark* 3.6. We have picked our grading convention as it is convenient for our present paper. It would be at least equally justified to define the spectrum $\mathrm{TJF}_m$ as $\Sigma^{-2m}\Gamma\mathcal{O}_{\mathcal{E}^{\mathrm{or}}}(me)$. This would better match the classical grading conventions for Jacobi forms (since $\mathbf{C} \otimes \pi_{2k}\Sigma^{-2m}\Gamma\mathcal{O}_{\mathcal{E}^{\mathrm{or}}}(me)$ are Jacobi forms of weight $k$ and index $\frac{m}{2}$) and might also be helpful to show an $E_\infty$-version of Theorem 3.3.

Let $P^m = \mathrm{cofib}(\Sigma\mathbf{C}P^{m-1} \to S^0)$ be the cofiber of the reduced $T$-transfer map, i.e. the suspension of the restriction of the map (A.9) to $\Sigma\mathbf{C}P^{m-1}$; we define $P^0$ as $S^0 \vee S^1$ by convention. By Proposition A.11, one can view $P^m$ as the spectrum "$\Sigma^2 \mathbf{C}P_{-1}^{m-1}$ without the 2-cell" for $m > 0$. We are now ready to prove Theorem 1.7.

**Theorem 3.7.** *There exists an equivalence of* TMF*-module spectra*

$$\mathrm{TJF}_m \simeq \mathrm{TMF} \wedge P^m.$$

*Under this equivalence, the map $a\colon L_m^{\mathrm{top}} \to L_{m+1}^{\mathrm{top}}$ corresponds to the map induced by the map $P^m \to P^{m+1}$ coming from the skeletal inclusion $\mathbf{C}P^{m-1} \to \mathbf{C}P^m$. In particular,*

$$\mathrm{TJF}_\infty = \mathrm{colim}_m \mathrm{TJF}_m = (a^{-1}\mathrm{TJF}_*)_0 \simeq \mathrm{TMF} \wedge P^\infty.$$

*Proof.* Recall from [GM23, Theorem 10.1] that the map

$$(3.8) \qquad (\mathrm{i}, \mathrm{tr})\colon \mathrm{TMF} \oplus \Sigma\,\mathrm{TMF} \to \mathrm{TMF}^T = (\mathrm{TMF}_T)^T$$

is an equivalence, where

- i is the unit map $\mathrm{TMF} \to \mathrm{TMF}^T$ (induced by the projection $p\colon \mathcal{E}^{\mathrm{or}} \to \mathcal{M}_{\mathrm{ell}}^{\mathrm{or}}$);
- $\mathrm{tr}\colon \Sigma\,\mathrm{TMF} \to \mathrm{TMF}^T$ is the degree-shifting transfer (cf. Appendix A).

Denoting by res the restriction $\mathrm{TMF}^T \to \mathrm{TMF}$ (induced by the zero section $e\colon \mathcal{M}_{\mathrm{ell}}^{\mathrm{or}} \to \mathcal{E}^{\mathrm{or}}$), note that $\mathrm{res} \circ \mathrm{i}$ is the identity.

On the other hand, $\mathrm{res} \circ \mathrm{tr}$ is the TMF-Hurewicz image of the map
$$S^1 \simeq \left(\Sigma_+^\infty T\right)^T \to (S^0)^T \to S^0,$$
where we use the Adams isomorphism Proposition A.5 in the first equivalence. For a general Lie group $G$ instead of $T$, this map always represents the homotopy class associated with the framed bordism class of $G$ in its invariant framing (cf. [Bau04]), thus we obtain

(3.9) $$\mathrm{res} \circ \mathrm{tr} \sim \eta$$

where $\eta \in \pi_1 \mathrm{TMF}$ is the Hurewicz image of the Hopf map. While not strictly necessary for this proof, we need this identification later.

Thus, we have a cofiber sequence
$$\Sigma \mathrm{TMF} \xrightarrow{\mathrm{tr}} \mathrm{TMF}^T \xrightarrow{\mathrm{res} - (0,\eta)} \mathrm{TMF},$$
where $(0,\eta)\colon \mathrm{TMF}^T \simeq \mathrm{TMF} \vee \Sigma \mathrm{TMF} \to \mathrm{TMF}$ is zero on the first component and $\eta$ on the second.

Consider the unit sphere $S(m\rho)$ in the representation $m\rho \cong \mathbf{C}^m$. This is a free $T$-space, and hence the Adams isomorphism in the form of Proposition A.5 gives an equivalence
$$\mathrm{TMF} \wedge \Sigma \mathbf{C}P_+^{m-1} \simeq \mathrm{TMF} \wedge \Sigma S(m\rho)_+/T \simeq (\mathrm{TMF}_T \wedge S(m\rho)_+)^T$$
By naturality of the Adams isomorphism, the inclusions $S(m\rho) \hookrightarrow S((m+1)\rho)$ correspond to the skeletal inclusions $\mathbf{C}P^{m-1} \hookrightarrow \mathbf{C}P^m$.

The cofiber sequence $S(m\rho)_+ \to S^0 \to S^{m\rho}$ gives rise to a cofiber sequence
$$\mathrm{TMF} \wedge \Sigma \mathbf{C}P_+^{m-1} \simeq (\mathrm{TMF}_T \wedge S(m\rho)_+)^T \to (\mathrm{TMF}_T)^T$$
$$\to (\mathrm{TMF}_T \wedge S^{m\rho})^T = \mathrm{TJF}_m\,.$$

This extends to a bigger diagram

$$\begin{array}{ccc}
\mathrm{TMF} \wedge \Sigma \mathbf{C}P_+^0 \simeq (\mathrm{TMF}_T \wedge S(\rho)_+)^T & \xrightarrow{\simeq} & \Sigma \mathrm{TMF} \\
\downarrow & & \downarrow \mathrm{tr} \\
\mathrm{TMF} \wedge \Sigma \mathbf{C}P_+^{m-1} \simeq (\mathrm{TMF}_T \wedge S(m\rho)_+)^T & \xrightarrow{\mathrm{tr}_m} (\mathrm{TMF}_T)^T & \longrightarrow \mathrm{TJF}_m \\
\downarrow & \downarrow \mathrm{res}-(0,\eta) & \downarrow = \\
\mathrm{TMF} \wedge \Sigma \mathbf{C}P^{m-1} & \xrightarrow{t} \mathrm{TMF} & \longrightarrow \mathrm{TJF}_m,
\end{array}$$

where we obtained the last row by taking cofibers of the vertical morphisms. Here, we have used that $S(\rho)_+ = T_+$ and the description of the degree-shifting transfer from Lemma A.8 for the commutativity of the upper square. Since the map $\Sigma \mathbf{C}P_+^{m-1} \to \Sigma \mathbf{C}P^{m-1}$ is split, we can compute $t$ as the composite
$$\mathrm{TMF} \wedge \Sigma \mathbf{C}P^{m-1} \to \mathrm{TMF} \wedge \Sigma \mathbf{C}P_+^{m-1} \xrightarrow{\mathrm{tr}_m} (\mathrm{TMF}_T)^T \xrightarrow{\mathrm{res}-(0,\eta)} \mathrm{TMF}\,.$$

We claim that the composite
$$\mathrm{TMF} \wedge \Sigma \mathbf{C}P^{m-1} \to \mathrm{TMF} \wedge \Sigma \mathbf{C}P_+^{m-1} \xrightarrow{\mathrm{tr}_m} (\mathrm{TMF}_T)^T \simeq \mathrm{TMF} \vee \Sigma \mathrm{TMF} \xrightarrow{\mathrm{pr}_2} \Sigma \mathrm{TMF}$$

is trivial. Accepting this claim for the moment, we deduce from Lemma A.12 that $t$ is TMF smashed with the reduced transfer $\Sigma \mathbf{C}P^{m-1} \to \Sigma \mathbf{C}P^{m-1}_+ \to S^0$ (since $(0, \eta)$ factors over $\mathrm{pr}_2$) and thus its cofiber agrees with $\mathrm{TMF} \wedge P^m$.

To show the claim, we apply Proposition 3.4 to the map $S(m\rho)_+ \to S^0$ and obtain that the dual of $\mathrm{tr}_m$ is equivalent to $\Sigma^{-1}$ of

$$(\mathrm{TMF}_T)^T \simeq \mathrm{TMF} \vee \Sigma \, \mathrm{TMF} \to \left(\mathrm{TMF}_T^{S(m\rho)_+}\right)^T \simeq \mathrm{TMF}^{\mathbf{C}P^{m-1}_+}.$$

This is a map of TMF-algebras. Thus, on the first summand the map is the unit map $\mathrm{TMF} \to \mathrm{TMF}^{\mathbf{C}P^{m-1}_+}$ induced by the map $q\colon \mathbf{C}P^{m-1}_+ \to S^0$ collapsing $\mathbf{C}P^{m-1}$ to a point. Dualizing back, $\mathrm{TMF} \wedge \mathbf{C}P^{m-1}_+ \xrightarrow{\mathrm{tr}_m} (\mathrm{TMF}_T)^T \xrightarrow{\mathrm{pr}_2} \Sigma \, \mathrm{TMF}$ must be $\pm \mathrm{TMF} \wedge q$.[2] Thus, the restriction of $\mathrm{pr}_2 \, \mathrm{tr}_m$ to $\mathrm{TMF} \wedge \mathbf{C}P^{m-1}$ is indeed trivial. $\square$

We will be making use of a slightly more general version of the transfer, along with stabilization and restriction maps, as follows.

The cofiber sequence of $T$-spaces

$$S(\rho)_+ \to S^0 \to S^\rho$$

considered before can be smashed with $S^{m\rho}$ to yield a cofiber sequence

$$S(\rho)_+ \wedge S^{m\rho} \to S^{m\rho} \to S^{(m+1)\rho}.$$

Since the left-hand term is still a free $T$-space, smashing with $\mathrm{TMF}_T$, taking $T$-fixed points, and applying the Adams isomorphism gives a cofiber sequence

(3.10) $\qquad \Sigma \, \mathrm{TMF} \wedge S^{2m} \xrightarrow{\mathrm{tr}} \mathrm{TJF}_m \xrightarrow{a} \mathrm{TJF}_{m+1} \xrightarrow{\mathrm{res}} \Sigma^{2m+2} \, \mathrm{TMF},$

where the map on the right is identified with the restriction map

$$(\mathrm{TMF}_T \wedge S^{(m+1)\rho})^T \to \mathrm{TMF}_T \wedge S^{(m+1)\rho}.$$

Assembling these triangles gives a sequence of distinguished triangles, infinite in both directions,

(3.11)
$$\cdots \xrightarrow{a} \mathrm{TJF}_{-1} \xrightarrow{a} \mathrm{TJF}_0 \xrightarrow{a} \mathrm{TJF}_1 \xrightarrow{a} \cdots$$
with res/tr maps down to $\Sigma^{-2}\,\mathrm{TMF}$, $\mathrm{TMF}$, $\Sigma^2\,\mathrm{TMF}$.

A struck-through arrows denote a map of degree $-1$. Truncating at $m \geq 0$ yields a similar sequence of distinguished triangles:

(3.12)
$$\mathrm{TJF}_0 \xrightarrow{a} \mathrm{TJF}_1 \xrightarrow{a} \mathrm{TJF}_2 \xrightarrow{a} \cdots$$
with $j_0$ (an equivalence), tr, $j_1 = \mathrm{res}$, tr, $j_2 = \mathrm{res}$ maps down to $\mathrm{TMF} \vee \Sigma\,\mathrm{TMF}$, $\Sigma^2\,\mathrm{TMF}$, $\Sigma^4\,\mathrm{TMF}$.

By definition of the truncation, we have $\mathrm{TJF}_0$ in the lower left corner. We will identify this with $\mathrm{TMF} \vee \Sigma\,\mathrm{TMF}$ via a map $j_0$ so that the map $\mathrm{res}\colon \mathrm{TJF}_0 \to \mathrm{TMF}$

---

[2]The sign comes from the fact that we did not specify the effect of the natural equivalence from Proposition 3.4 in the case of $X = S^0$, which takes the form

$$\mathrm{TMF}[1] \oplus \mathrm{TMF} \simeq \Sigma D_{\mathrm{TMF}}(\mathrm{TMF} \oplus \mathrm{TMF}[1]) \simeq \mathrm{TMF} \oplus \mathrm{TMF}[1].$$

Using that $\pi_{-1}\,\mathrm{TMF} = 0$ and that the units in $\pi_0\,\mathrm{TMF} = \mathbf{Z}[j]$ are $\pm 1$, the corresponding matrix must be of the form $\begin{pmatrix} ? & \pm 1 \\ \pm 1 & 0 \end{pmatrix}$, which suffices for our purposes.

becomes $\mathrm{pr}_1$ under this identification. Concretely, let $1 \in \pi_0 \mathrm{TJF}_0$ be the multiplicative unit, $z \in \pi_2 \Sigma^2 \mathrm{TMF}$ the double suspension of the unit and $\tau = \mathrm{tr}(z) \in \pi_1 \mathrm{TJF}_0$ its image under the transfer map $\mathrm{tr} \colon \Sigma \mathrm{TMF} \to \mathrm{TJF}_0$. By (3.8), 1 and $\tau$ generate $\pi_* \mathrm{TJF}$ freely as a TMF-module; denoting the standard basis of $\mathrm{TMF} \vee \Sigma \mathrm{TMF}$ also by 1 and $\tau$, we set $j_0(1) = 1$ and $j_0(\tau) = \eta + \tau$. The relation $\mathrm{res}(\mathrm{tr}(z)) = \eta$ shows that indeed $\mathrm{pr}_1 j_0 = \mathrm{res}$.

Since $\bigvee_m \mathrm{TJF}_m$ is a filtered graded ring spectrum (with $j_0$ inducing the multiplication on $\mathrm{TMF} \vee \Sigma \mathrm{TMF}$), both (3.11) and (3.12) give rise to multiplicative exact couples. In particular, the following relations hold for $x \in \pi_* \Sigma^{2m} \mathrm{TMF}$, $y \in \pi_* \mathrm{TJF}_n$:

$$\text{(3.13)} \qquad \mathrm{tr}(x \operatorname{res}(y)) = \mathrm{tr}(x) y \in \pi_* \mathrm{TJF}_{m+n}.$$

We can now identify the ring $\pi_* \mathrm{TJF}_0$ as $\pi_* \mathrm{TMF}[\tau]/(\tau^2 - \tau \eta)$:

*Proof of Prop. 1.10.* For $m = 0$, the restriction map res in (3.10) is trivial and $a$ is split by the unit $\mathrm{TMF} \to \mathrm{TJF}_0$, as shown in [GM23, Theorem 10.1]. Then (3.13) and (3.9) imply

$$\tau^2 = \mathrm{tr}(1)\tau = \mathrm{tr}(1 \operatorname{res}(\tau)) = \mathrm{tr}(\operatorname{res}(\tau)) = \eta \tau. \qquad \square$$

In homotopy, diagram (3.11) gives rise to a spectral sequence with

$$E^1_{**} = \pi_* \mathrm{TMF}[z^{\pm 1}]$$

and convering to the homotopy of the spectrum $\mathrm{TMF}^{\mathbf{C}P^\infty}$ considered in [AFG08]. We are more interested in the other spectral sequence, converging to the homotopy of $\mathrm{TJF}_\infty$:

**Proposition 3.14.** *The $E^1$-term of the homotopy spectral sequence arising from (3.12) is given by*

$$E^1_{**} = \pi_* \mathrm{TMF}[z, \tau]/(\tau^2 - \tau \eta, z\tau - z\eta)$$

*Proof.* Let us temporarily denote by $\bar{\tau}$ the element

$$j_0(\tau) = (\eta, 1) \in \pi_1 \mathrm{TMF} \vee \Sigma \mathrm{TMF}.$$

Since $j_0$ is by definition multiplicative, the equation $\bar{\tau}^2 = \bar{\tau} \eta$ holds. We want to show $z\bar{\tau} = z\eta \in \pi_3 \Sigma^2 \mathrm{TMF}$.

Since $j_1 \colon \mathrm{TJF}_1 \simeq \mathrm{TMF} \to \Sigma^2 \mathrm{TMF}$ is null-homotopic by the splitting, the map $\mathrm{tr} \colon \Sigma \mathrm{TMF} \to \mathrm{TJF}_0$ is injective in homotopy and it suffices to show $\mathrm{tr}(z\bar{\tau}) = \mathrm{tr}(z\eta)$. But

$$\mathrm{tr}(z\bar{\tau}) = \mathrm{tr}(z j_0(\tau)) = \mathrm{tr}(z)\tau = \tau^2 = \eta \tau = \mathrm{tr}(z\eta).$$

Thus $z\bar{\tau} = z\eta$. The claim follows by abusing notation and writing $\tau$ for $\bar{\tau}$. $\square$

Using the unit map $i \colon \mathrm{TMF} \to \mathrm{TJF}_0$, we can define a modified ring spectrum of topological Jacobi forms whose $m$-th graded piece is

$$\overline{\mathrm{TJF}}_m = \mathrm{TMF} \times_{\mathrm{TJF}_0} \mathrm{TJF}_m.$$

This is an additive (but not multiplicative) factor of the graded ring spectrum TJF. The corresponding modified sequence of distinguished triangles has the form

(3.15)
$$\overline{\mathrm{TJF}}_0 \xrightarrow[\simeq]{a} \overline{\mathrm{TJF}}_1 \xrightarrow{a} \overline{\mathrm{TJF}}_2 \xrightarrow{q} \cdots$$

with maps $j_0 \simeq$ down to TMF, then to $0$, then tr from $\Sigma^4 \mathrm{TMF}$, and $j_2 = \mathrm{res}$.

Since the standard inclusion (3.15) ↪ (3.12) is multiplicative, the associated spectral sequence (converging to $\mathrm{TMF}_\infty$ as well) has

$$E^1_{**} = \pi_* \operatorname{TMF}[x,y]/(x^3 - y^2),$$

where $x$ maps to $z^2$ and $y$ to $z^3$. Identifying $\mathrm{TJF}_m$ with $\mathrm{TMF} \wedge P^m$ and noting that the map $a$ is induced by the inclusion $P^m \to P^{m+1}$ we can also see that the filtration (3.15) exactly corresponds to the skeletal filtration of $P^\infty$ and the spectral sequence above is isomorphic to the Atiyah-Hirzebruch spectral sequence. This will allow us to easily derive differentials.

## 4. Computations of derived Jacobi forms

In this purely algebraic section, we will compute the cohomology of the (stacky, but not spectral) universal elliptic curve $\mathcal{E}$ over $\mathcal{M}_{\mathrm{ell}}$ with coefficients in the Looijenga line bundle $L_{2n,m} = p^*\omega^n \otimes \mathcal{O}(me)$ for all $m$ and $n$, including the case $m = \infty$, which we agree to mean $L_{2n,\infty} = p^*\omega^n \otimes \operatorname{colim}_m \mathcal{O}(me)$. In fact, we will instead consider the larger stack $\mathcal{M}_{\mathrm{Weier}}$ of curves with a Weierstrass parametrization; the stack $\mathcal{M}_{\mathrm{ell}}$ is the substack of $\mathcal{M}_{\mathrm{Weier}}$ where the discriminant $\Delta$ is invertible: $\mathcal{M}_{\mathrm{ell}} \cong \Delta^{-1}\mathcal{M}_{\mathrm{Weier}}$. Similarly, the universal elliptic curve $\mathcal{E}$ is replaced with the universal Weierstrass curve $\mathcal{W} \to \mathcal{M}_{\mathrm{Weier}}$, which is singular over the complement $\mathcal{M}_{\mathrm{Weier}} - \mathcal{M}_{\mathrm{ell}}$. While $\mathcal{W}$ is a stack, it is also a relative projective scheme in the sense that for any affine $E\colon \operatorname{Spec} R \to \mathcal{M}_{\mathrm{Weier}}$, the basechange $\mathcal{W}_E = E^*\mathcal{W}$ is projective. We will also consider the relative affine scheme $\mathcal{W}_0 = \mathcal{W} - \{e\}$ obtained by removing the zero section (the point at infinity). Note that sections of $L_{2n,\infty}$ over $f\colon U \to \mathcal{W}$ are precisely sections of $p^*\omega^n$ over $f^{-1}(\mathcal{W}_0)$.

The Looijenga line bundle $L_{n,m}$ extends to $\mathcal{W}$ because $\omega$ extends to $\mathcal{M}_{\mathrm{Weier}}$ (it has a coordinate given by $\omega = \frac{dx}{2y+a_1x+a_3}$ for a curve in Weierstrass form), and $\mathcal{O}(me)$ still makes sense on $\mathcal{W}$. Neither $\mathcal{M}_{\mathrm{Weier}}$ nor $\mathcal{W}$ are Deligne-Mumford stacks, but they are Artin stacks. They have indeed flat affine covers by affine schemes, and thus the cohomology $H^*(\mathcal{W}, L_{n,m})$ can be computed as the cohomology of Hopf algebroids. We briefly recall this approach.

Let $\mathcal{M}$ be an algebraic stack with a flat affine cover $\operatorname{Spec} A \xrightarrow{p} \mathcal{M}$. Then $\operatorname{Spec} A \times_\mathcal{M} \operatorname{Spec} A = \operatorname{Spec} \Gamma$ is affine as well, and $(A, \Gamma)$ is a Hopf algebroid (it represents an affine groupoid scheme). Moreover, if $\mathcal{F}$ is a quasicoherent $\mathcal{O}_\mathcal{M}$-module sheaf on $\mathcal{M}$ then $p^*\mathcal{F}$ is the sheaf associated to an $A$-module $M$, which more strongly is a (left, say) $(A, \Gamma)$-comodule: there is a coaction $M \xrightarrow{\psi} \Gamma \otimes_A M$. Flatness ensures that there are isomorphisms

$$H^*(\mathcal{M}, \mathcal{F}) \cong \operatorname{Ext}^*_{(A,\Gamma)}(A, M) =: H^*(A, \Gamma; M),$$

where the Ext group is in the abelian category of $(A, \Gamma)$-comodules. If, moreover, $\mathcal{F}$ is a quasicoherent commutative $\mathcal{O}_\mathcal{M}$-algebra sheaf, then the comodule $M$ is a comodule algebra, i.e. it has a unit $u\colon A \to M$ and a product $M \otimes_A M \to M$ which are maps of comodules, where $M \otimes_A M$ becomes a comodule by

$$M \otimes_A M \xrightarrow{\psi \otimes \psi} \Gamma \otimes_A M \otimes_A \Gamma \otimes_A M \to \Gamma \otimes_A \Gamma \otimes_A M \otimes_A M \xrightarrow{\mu \otimes \mathrm{id} \otimes \mathrm{id}} \Gamma \otimes_A M \otimes_A M,$$

and the cohomology isomorphism above is multiplicative. In fact, in this case, we can construct another Hopf algebroid

$$H_M = (M, \Gamma \otimes_A M)$$

with the left and right units given by $\eta_R \otimes \mathrm{id}$ and $\psi$, respectively, the augmentation $\Gamma \otimes_A M \to M$ by $\epsilon \otimes \mathrm{id}$, and the comultiplication by

$$\Gamma \otimes_A M \xrightarrow{\Delta \otimes \mathrm{id}} \Gamma \otimes_A \Gamma \otimes_A M \cong (\Gamma \otimes_A M) \otimes_M (\Gamma \otimes_A M).$$

This Hopf algebroid agrees with the simple base-change Hopf algebroid $(B, \Gamma \otimes_A B)$ for a $B$-algebra $A$ in the case where $B = M$ is a trivial comodule, i.e. $\psi = \eta_L \otimes \mathrm{id}$. We have the following very special case of the change-of-rings isomorphism of [Rav04, A1.3.12]:

**Lemma 4.1.** *Let $H = (A, \Gamma)$ be a Hopf algebroid and $M$ a commutative $H$-comodule algebra such that $M$ is flat as an $A$-module. Let $H_M$ denote the Hopf algebroid $(M, \Gamma \otimes_A M)$ constructed above. Then for any $H_M$-comodule $C$, there is a natural isomorphism*

$$\mathrm{Ext}^*_H(A, C) \cong \mathrm{Ext}^*_{H_M}(M, C). \qquad \square$$

*Convention* 4.2. In the following, all our Hopf algebroids are understood to be *graded*. A graded Hopf algebroid $(A, \Gamma)$ induces the structure of an ungraded Hopf algebroid on the pair $(A, \Gamma[u^{\pm 1}])$, where the powers of $u$ track the grading. By definition, the stack associated to a graded Hopf algebroid is the stack associated to the associated ungraded Hopf algebroid. See [MO20, Section 4] for details.

Here is a short summary of the various stack cohomologies and Hopf algebroids used in this section.

- The bigraded ring of *derived modular forms*

$$\mathrm{dmf}_{2t-s,s} = H^s(\mathcal{M}_{\mathrm{Weier}}, \omega^t)$$

as originally computed by Hopkins and Mahowald (cf. [Bau08]). The subring $\mathrm{dmf}_{2t,0}$ is the ring $\mathrm{mf}_{2t}$ of classical modular forms of weight $t$ over $\mathbf{Z}$. The ring $\mathrm{dmf}_{**}$ forms the $E_2$-term of the descent spectral sequence $E_2^{s,t} = \mathrm{dmf}_{2t-s,s} \Rightarrow \pi_{2t-s} \mathrm{tmf}$. The stack $\mathcal{M}_{\mathrm{Weier}}$ has a presentation by the Hopf algebroid $H = (A, \Gamma)$ with $A = \mathbf{Z}[a_1, a_2, a_3, a_4, a_6]$ and $\Gamma = A[r, s, t]$. The cover $\mathrm{Spec}\, A \to \mathcal{M}_{\mathrm{Weier}}$ classifies the Weierstrass curve given by the affine equation

(4.3) $\quad E(a_1, a_2, a_3, a_4, a_6, x, y) = y^2 + a_1 xy + a_3 - (x^3 + a_2 x^2 + a_4 x + a_6) = 0.$

- The bigraded filtered $\mathrm{dmf}_{**}$-algebra of *derived, stable, connective Jacobi forms*

$$\mathrm{djF}_{2t-s,s,\infty} = H^s(\mathcal{W}_0, p^*(\omega^t)).$$

As a stack, $\mathcal{W}_0$ has a presentation by the Hopf algebroid

$$H_{\mathcal{W},\infty} = (A_\mathcal{W}, A_\mathcal{W} \otimes_A \Gamma) = (A_\mathcal{W}, A_\mathcal{W}[r, s, t])$$

with

$$A_\mathcal{W} = A[x, y]/E(a_1, \ldots, a_6, x, y).$$

The structure maps of $H_{\mathcal{W},\infty}$ restrict to the structure maps of $H$ under the inclusion $A \to A_\mathcal{W}$, and

$$\eta_L(x) = x;\ \eta_L(y) = y;\ \eta_R(x) = x - r;\ \eta_R(y) = y - sx + sr - t.$$

Note that $\eta_R$ describes the inverse transformation of $(x, y) \mapsto (x + r, y + sx + t)$.

The category of quasicoherent $\mathcal{O}_{\mathcal{W}_0}$-module sheaves is trivially equivalent to the category of quasicoherent $\mathcal{O}_{\mathcal{W}}(\infty e)$-module sheaves, which means that as dmf-algebras,
$$\mathrm{djF}_{t-s,s,\infty} \cong H^s(\mathcal{W}, L_{t,\infty})$$
Since $L_{t,\infty} = \operatorname{colim}_m L_{t,m}$ by definition, the canonical maps
$$i_m \colon H^s(\mathcal{W}, L_{t,m}) \to \mathrm{djF}_{t-s,s,\infty}$$
induce an exhaustive, increasing filtration $F_m \mathrm{djF}_{t-s,s,\infty} = \operatorname{im}(i_m)$ by dmf-submodules. Moreover, this filtration is multiplicative in the sense that $(F_m \mathrm{djF}_{*,*,\infty})(F_{m'} \mathrm{djF}_{*,*,\infty}) \subseteq F_{m+m'} \mathrm{djF}_{*,*,\infty}$. In terms of the Hopf algebroid $H_{\mathcal{W},\infty}$, this filtration is realized as a filtration
$$F_* H_{\mathcal{W},\infty} = (F_* A_{\mathcal{W}}, F_* A_{\mathcal{W}} \otimes_A \Gamma)$$
where $x \in F_2 A_{\mathcal{W}} - F_1 A_{\mathcal{W}}$, $y \in F_3 A_{\mathcal{W}} - F_2 A_{\mathcal{W}}$ and $a_i \in F_0 A_{\mathcal{W}}$ for all $i$ because the coordinate function $x$ has a double pole at $\infty$, while the coordinate function $y$ has a triple pole.

- The trigraded dmf$_{**}$-algebra of (unstable) *derived, connective Jacobi forms*
$$\mathrm{djF}_{t-s,s,m} = H^s(\mathcal{W}, L_{t,m}),$$
where $L_{t,m}$ is the Looijenga line bundle, which over $\mathcal{M}_{\mathrm{ell}}$ is realized as the homotopy of the sheaf $\mathcal{O}_{\mathcal{E}^{\mathrm{or}}}(me)$ of Thm. 3.3. Thus djF is indeed an algebra. For the projection $p \colon \mathcal{W} \to \mathcal{M}_{\mathrm{Weier}}$ we have that $R^i p_* \mathcal{O}_{\mathcal{W}}(me) = 0$ for $i > 0$ and $m > 0$, and hence
$$\mathrm{djF}_{t-s,s,m} \cong \begin{cases} H^s(\mathcal{M}_{\mathrm{Weier}}, p_* L_{m,t}); & m > 0 \\ \left( H^*(\mathcal{M}_{\mathrm{Weier}}, \omega^{\frac{t}{2}})[\tau]/(\tau^2 - \tau h_1) \right)_{s,t}; & m = 0, \end{cases}$$
the $m = 0$ case being part of Prop. 5.5 proved below.

We define the following dmf$_{**}$-subalgebra of $\mathrm{djF}_{***}$:
$$\overline{\mathrm{djF}}_{t-s,s,m} = H^s(\mathcal{M}_{\mathrm{Weier}}, p_* L_{m,t}),$$
which differs from djF only in index $m = 0$. As an $(A, \Gamma)$-comodule algebra, $p_* L_{*,*}$ is given by
$$A_{\mathcal{W},*} = A[a, x, y]/E'(a, a_1, \ldots, a_6, x, y),$$
where
$$(4.4) \quad E'(a, a_1, \ldots, a_6, x, y) = y^2 + a_1 a x y + a_3 a^3 y - (x^3 + a_2 a^2 x^2 + a_4 a^4 x + a_6 a^6)$$
and the bidegrees of the generators are $|a_i| = (2i, 0)$, $|a| = (0, 1)$, $|x| = (4, 2)$, $|y| = (6, 3)$. The powers of $a$ are forced by the fact that $E'$ has to reduce to $E$ (4.3) when setting $a = 1$ and homogeneity. Multiplication by $a$ represents the canonical inclusion $L_{m,n} \to L_{m+1,n}$. The coaction $A_{\mathcal{W},*} \to \Gamma \otimes_A A_{\mathcal{W},*} = A_{\mathcal{W},*}[r,s,t]$ is given by $\psi(a_i) = \eta_R(a_i)$, $\psi(a) = 1 \otimes a$, $\psi(x) = x - a^2 r$, and $\psi(y) = y - asx + a^2(sr - t)$. By Lemma 4.1, $\overline{\mathrm{djF}}_{*,*,*}$ is thus isomorphic to the cohomology of the Hopf algebroid $H_{\mathcal{W},*} = (A_{\mathcal{W},*}, \Gamma \otimes_A A_{\mathcal{W},*})$.

The map $A_{\mathcal{W},*} \to A_{\mathcal{W}}$, sending $a$ to 1, induces an isomorphism of Hopf algebroids
$$\left(a^{-1} H_{\mathcal{W},*}\right)_{m,0} \to (H_{\mathcal{W},\infty})_m,$$

which allows us to compare the cohomology of $H_{\mathcal{W},*}$ to the "stable" cohomology, i.e. that of $H_{\mathcal{W}}$, which is significantly easier to compute:

$$(a^{-1}\,\mathrm{djF}_{*,*,*})_{n,s,0} \cong (a^{-1}\overline{\mathrm{djF}}_{*,*,*})_{n,s,0} \cong \mathrm{djF}_{n,s,\infty}.$$

The first isomorphism holds because $\tau a = 0$.

Are there tmf-algebra spectra $\mathrm{tjF}_\infty$ and $\mathrm{tjF}_*$ whose homotopy is the abutment of descent spectral sequences with $E_2$-terms $\mathrm{djF}_{*,*,\infty}$ and $\mathrm{djF}_{*,*,*}$, respectively? This remains conjectural for the time being. We do know two things:

(1) The answers are yes if one relaxes the requirement to tmf-module spectra. In that case, $\mathrm{tjF}_m = \mathrm{tmf} \wedge P^m$ for $0 \le m \le \infty$.
(2) The answers are also positive when restricting to $\mathcal{M}_{\mathrm{ell}} \subseteq \mathcal{M}_{\mathrm{Weier}}$ and $\mathcal{E} \subseteq \mathcal{W}$ by the construction of $\mathrm{TJF}_m$ in Section 3.2.

In this paper, we will compute the following:

(1) $(\mathrm{djF}_{*,*,\infty})_{(p)}$ and $(\mathrm{djF}_{*,*,*})_{(p)}$ for $p \ge 3$;
(2) $(\mathrm{djF}_{*,*,\infty})_{(2)}$, including its filtration.
(3) $\pi_*(\mathrm{tjF}_\infty)_{(p)}$ for all $p$ and $\pi_*(\mathrm{tjF}_*)_{(p)}$ for $p \ge 3$ with the definition $\mathrm{tjF}_m = \mathrm{tmf} \wedge P^m$ as above. One obtains the corresponding results for $\mathrm{TJF}_m$ by inverting $\Delta$.

The complete results are summarized in the following theorems. These are somewhat stronger than Thms. 1.9 and 1.11 in that our computations here deal with weak derived Jacobi forms rather than weakly holomorphic ones.

The first computation is the connective version of Thm. 1.9 away from 2:

**Theorem 4.5.**
$$(\mathrm{djF}_{*,s,\infty})[\tfrac{1}{2}] \cong \begin{cases} \mathbf{Z}[\tfrac{1}{2}, b_2, b_3, b_4]; & s = 0 \\ 0; & s > 0 \end{cases}$$

with $|b_i| = 2i$. The algebra structure over $(\mathrm{dmf}_*)$ factors through

$$(\mathrm{mf}_*)[\tfrac{1}{2}] = \mathbf{Z}[\tfrac{1}{2}, c_4, c_6, \Delta]/(c_4^3 - c_6^2 - 1728\Delta)$$

and is given by

$$c_4 \mapsto b_2^2 - 24b_4$$
$$c_6 \mapsto -b_2^3 + 36b_2 b_4 - 216 b_3^2$$
$$\Delta \mapsto \frac{1}{4}b_2^3 b_3^2 - 27 b_3^4 + 9 b_2 b_3^2 b_4 + \frac{1}{4} b_2^2 b_4^2 - 8 b_4^3.$$

The associated graded of the filtration is given by

$$\mathrm{gr}_* \mathrm{djF}_{*,*,\infty}[\tfrac{1}{2}] = \mathrm{mf}_*[\tfrac{1}{2}, b_2, b_3, b_4]/I,$$

where

$$I = (b_2^2 - 24 b_4, -b_2^3 + 36 b_2 b_4 - 216 b_3^2, \frac{1}{4} b_2^3 b_3^2 - 27 b_3^4 + 9 b_2 b_3^2 b_4 + \frac{1}{4} b_2^2 b_4^2 - 8 b_4^3)$$

with $b_i$ in the $m = 2i$th filtration quotient.

Moreover, $\pi_*(\mathrm{tjF}_\infty)[\tfrac{1}{2}] \cong (\mathrm{djF}_{*,0,\infty})[\tfrac{1}{2}]$.

Next we verify Thm. 1.11 in the connective setting:

**Theorem 4.6.**
$$(\mathrm{djF}_{*,*,*})[\tfrac{1}{2}] \cong \mathbf{Z}[\tfrac{1}{2}, a, b_2, b_3, b_4, c_4, c_6, \Delta, \tau, \alpha, \beta, \gamma]/I,$$
*with the following tridegrees:*
$$|a| = (0,0,1); \; |b_i| = (2i, 0, i);$$
$$|c_4| = (8,0,0); \; |c_6| = (12,0,0); \; |\Delta| = (24,0,0);$$
$$|\tau = (1,1,0); \; |\alpha| = (3,1,0); \; |\beta| = (10,2,0); \; |\gamma| = (7,1,2),$$
*and where the ideal $I$ is*
$$I = (1728\Delta - c_6^2 + c_4^3, c_4\alpha, c_4\beta, c_4\gamma, c_6\alpha, c_6\beta, c_6\gamma,$$
$$a^4 c_4 - b_2^2 + 24 b_4, a^6 c_6 + b_2^3 - 36 b_2 b_4 + 216 b_3^2,$$
$$a^{12}\Delta - \frac{1}{4} b_2^3 b_3^2 + 27 b_3^4 - 9 b_2 b_3^2 b_4 - \frac{1}{4} b_2^2 b_4^2 + 8 b_4^3,$$
$$\tau^2, \tau a, \tau b_2, \tau b_3 - 2a\gamma, \tau b_4, \tau\gamma$$
$$\alpha^2, 3\alpha, 3\beta, b_2\alpha, b_4\alpha, b_2\beta, b_4\beta, 3\gamma, a^2\gamma, \gamma^2, a^2\beta - \gamma\alpha, b_2\gamma, b_4\gamma, a^2\alpha).$$

*The $(\mathrm{dmf}_{**})[\tfrac{1}{2}]$-algebra structure is given by the fact that $\alpha \in \mathrm{dmf}_{3,1}$ and $\beta \in \mathrm{dmf}_{10,2}$ and the modular forms $c_4$, $c_6$, $\Delta$ map to the classes of the same name in* djF.

Finally, we complete the proof of Thm. 1.9 with the statement localized at $p = 2$:

**Theorem 4.7.**
$$(\mathrm{djF}_{*,*,\infty})_{(2)} \cong \mathbf{Z}_{(2)}[b_2, b_3, b_4, b_8, h_1]/(2h_1, b_3 h_1, 4b_8 + b_4^2 - b_2 b_3^2)$$
*with*
$$|b_i| = (2i, 0); \; h_1 = (1,1).$$
*The $(\mathrm{dmf}_{**})$-module structure is given by $c_4 \mapsto b_2^2 - 24 b_4$, $c_6 \mapsto -b_2^3 + 36 b_2 b_4 - 216 b_3^2$, $\Delta \mapsto -b^2 b_8 - 8 b_4^3 - 27 b_3^4 + 9 b_2 b_3^2 b_4$, $h_1 \mapsto h_1$, $\delta \mapsto b_2 h_1$, and $h_2, \epsilon, \kappa, \bar\kappa \mapsto 0$.*

## 5. Computations at odd primes

When localized at a prime $p$, the stack $\mathcal{M}_{\mathrm{Weier}}$ has a simpler presentation (cf. [Bau08]), and correspondingly, so does $\mathcal{W}$. In terms of Hopf algebroids, if $(A, \Gamma) \to (A', \Gamma')$ is an equivalence of Hopf algebroids then so is $(B, B \otimes_A \Gamma) \to (A' \otimes_A B, A' \otimes_A \Gamma')$.

Localizing away from 2, we can replace $y$ in the affine equation (4.3) by $\tfrac{1}{2}(y - a_1 x - a_3)$ to transform it into the form
$$y^2 = 4x^3 + b_2 x^2 + 2b_4 x + b_6.$$
The Hopf algebroid $H_{\mathcal{W},\infty}[\tfrac{1}{2}]$ becomes equivalent with the Hopf algebroid $\tilde{H}_\infty = (\tilde{A}_\infty, \tilde{\Gamma}_\infty)$ with $\tilde{A}_\infty = \mathbf{Z}[\tfrac{1}{2}, b_2, b_4, b_6, x, y]/(y^2 - 4x^3 - b_2 x^2 - 2b_4 x - b_6)$ and $\tilde{\Gamma}_\infty = \tilde{A}[r]$ (cf. [Bau08, Lemma 4.1]). Similarly, the replacement of $y$ by $\tfrac{1}{2}(y - a_1 a x - a_3 a^3)$ transforms (4.4) into
$$y^2 = 4x^3 + b_2 a^2 x^2 + 2 b_4 a^4 x + b_6 a^6$$
and $H_{\mathcal{W},*}$ becomes equivalent with $\tilde{H}_* = (\tilde{A}_*, \tilde{\Gamma}_*)$ with
$$\tilde{A}_* = \mathbf{Z}[\tfrac{1}{2}, a, b_2, b_4, b_6, x, y]/(y^2 - 4x^3 - b_2 a^2 x^2 - 2 b_4 a^4 x - b_6 a^6) \text{ and } \tilde{\Gamma}_* = \tilde{A}_*[r].$$

In each case, the right units of the Hopf algebroids are given by (set $a = 1$ for the case $A_\infty$):

$$\eta_R(b_2) = b_2 + 12r$$
$$\eta_R(b_4) = b_4 + b_2 r + 6r^2$$
$$\eta_R(b_6) = b_6 + 2b_4 r + b_2 r^2 + 4r^3$$
$$\eta_R(a) = a$$
$$\eta_R(x) = x - a^2 r$$
$$\eta_R(y) = y.$$

We make repeated use of the following change-of-cover theorem by Hopkins [HS99, Hov02]:

**Theorem 5.1.** *Let $(A, \Gamma)$ be a Hopf algebroid, $A \to A'$ a morphism of rings, and $\Gamma' = A' \otimes_A \Gamma \otimes_A A'$ be the base change of $\Gamma$ to a Hopf algebroid over $A'$. If there exists a ring $R$ and a morphism $A' \otimes_A \Gamma \to R$ such that the composite*

$$A \xrightarrow{1 \otimes \eta_R} A' \otimes_A \Gamma \to R$$

*is faithfully flat, then it induces an equivalence of comodule categories, and in particular an isomorphism*

$$H^*(A, \Gamma) \xrightarrow{\cong} H^*(A', \Gamma') \qquad \square$$

In the stable case, we find an easy application of this theorem by dividing out $x$.

**Lemma 5.2.** *The map $\tilde{A} \to A' = \tilde{A}/(x)$ induces an isomorphism*

$$\mathrm{djF}_{*,s,\infty}[\tfrac{1}{2}]) \cong \begin{cases} \mathbf{Z}[\tfrac{1}{2}, b_2, b_4, y]; & s = 0 \\ 0; & s > 0, \end{cases}$$

*Proof.* Let us first show that Thm. 5.1 is applicable for

$$R = A' \otimes_{\tilde{A}} \Gamma,$$

i.e. that

$$g = (1 \otimes \eta_R) \colon \tilde{A} \to A' \otimes_{\tilde{A}} \tilde{A}[r] = A'[r]$$

is faithfully flat (and actually an isomorphism). This follows from the fact that $\eta_R(x) = x - r$ and thus $g(x) = -r$ is monic up to a unit.

Now by definition,

$$A' = \mathbf{Z}[\tfrac{1}{2}, b_2, b_4, b_6, y]/(y^2 - b_6)$$

and

$$\Gamma' = \tilde{A}_\mathcal{W}[r]/(\eta_L(x), \eta_R(x)) = \tilde{A}_\mathcal{W}[r]/(x, r) = A'. \qquad \square$$

*Proof of Thm. 4.5.* Setting $b_3 = y$, the lemma above shows everything except for the claim about the filtration. The formulas for $c_4$, $c_6$ and $\Delta$ are the standard expressions computed, e.g., in [Sil09, Chapter III.1]. Unfortunately, our way of manipulating the Hopf algebroid $H_{\mathcal{W},\infty}$ was incompatible with the filtration, so we postpone the claim about the associated graded of the filtration as Cor. 5.8, a consequence of Thm. 4.6. $\qquad \square$

The computation required for Thm. 4.6 is a bit more complicated. The above computation does not immediately generalize because dividing out by $x$ no longer produces an equivalent Hopf algebroid.

We will first construct variants of the algebraic Atiyah-Hirzebruch spectral sequence for $H^*(\mathcal{W}, \mathcal{O}_\mathcal{W}(me))$. These are algebraic versions of (3.11), (3.12), and (3.15).

Consider the multiplicative filtration of $\mathcal{O}_\mathcal{W}$-module sheaves

(5.3) $$\cdots \xrightarrow{a} \mathcal{O}_\mathcal{W}(-e) \xrightarrow{a} \mathcal{O}_\mathcal{W} \xrightarrow{a} \mathcal{O}_\mathcal{W}(e) \xrightarrow{a} \cdots \hookrightarrow \mathcal{O}_\mathcal{W}(\infty e).$$

Note that for each $m \in \mathbf{Z}$, we have short exact sequences of sheaves of $\mathcal{O}_\mathcal{W}$-modules

$$0 \to \mathcal{O}_\mathcal{W}((m-1)e) \xrightarrow{a} \mathcal{O}_\mathcal{W}(me) \to e_*(\omega^{-m}) \to 0$$

which, when restricted to $\mathcal{E} \subset \mathcal{W}$ and $m = -1$, coincides with the homotopy sheaves of (3.2).

We can tensor with the vector bundle $p^*(\omega^t)$ to obtain a short exact sequence

$$0 \to L_{2t,m-1} \xrightarrow{a} L_{2t,m} \to e_*(\omega^{-m}) \otimes p^*(\omega^t) \to 0.$$

By the projection formula,

$$e_*(\omega^{-m}) \otimes p^*(\omega^t) \cong e_*(\omega^{-m} \otimes e^*p^*\omega^t) = e_*(\omega^{t-m}).$$

Since $e$ is a closed embedding, $e_*$ is exact and

$$H^s(\mathcal{W}, e_*(\omega^{t-m})) \cong H^s(\mathcal{M}_{\text{Weier}}, \omega^{t-m}) \cong \mathrm{dmf}_{2(t-m)-s,s}.$$

Upon applying $H^*(\mathcal{W}, -)$ and setting $n = 2t - s$, we thus obtain an unrolled exact couple

(5.4)
$$\begin{array}{ccccccc}
\cdots \xrightarrow{a} \mathrm{djF}_{n,s,-1} & \xrightarrow{a} & \mathrm{djF}_{n,s,0} & \xrightarrow{a} & \mathrm{djF}_{n,s,1} & \xrightarrow{a} & \cdots \\
& \searrow^{\text{res}} \nwarrow^{\text{tr}} & & \searrow^{\text{res}} \nwarrow^{\text{tr}} & & \searrow^{\text{res}} \nwarrow^{\text{tr}} & \\
\mathrm{dmf}_{n-2,s} & & \mathrm{dmf}_{n,s} & & \mathrm{dmf}_{n+2,s} & & \mathrm{dmf}_{n+4,s}
\end{array}$$

with connecting homomorphism $\mathrm{tr}\colon \mathrm{dmf}_{n-2m,s} \to \mathrm{djF}_{n-1,s+1,m-1}$ and thus a spectral sequence with

$$E^1_{***} = \mathrm{dmf}_{**}[z^{\pm 1}]$$

where the grading is such that $E^1_{s,2t,0} = \mathrm{dmf}_{2t-s,s}$ and the tridegree of $z$ is $(s, t, u) = (0, 0, 2)$.

The spectral sequence converges to a $T$-Tate variant of $\mathrm{djF}_{**\infty}$.

Upon inverting $\Delta$, both the $E^1$ and $D^1$-terms of (5.4) are themselves $E^2$-terms of descent spectral sequences converging to $\pi_*\Sigma^{2m}\mathrm{TMF}$ and $\pi_*\mathrm{TJF}_m$, respectively. Morally, (5.4) can be thought of as an exact couple of spectral sequences, but this is not straightforward since spectral sequences do not form an abelian category. We will leave it at that but note that the maps $a$, res, and tr in (5.4) converge to the maps of the same names between $\mathrm{TJF}_{m-1}$, $\mathrm{TJF}_m$, and $\Sigma^{2m}\mathrm{TMF}$.

We can truncate the filtration (5.3), considering only the terms $m \geq 0$; here we replace $\mathrm{dmf}_{2t-s,s}$ by $\mathrm{djF}_{2t-s,s,0}$ as $E^1_{s,2t,0}$ and posit that the map

$$\mathrm{dmf}_{*,*} \oplus \mathrm{dmf}_{*-1,*,-1} \cong D^1_{*,*,0} \xrightarrow{\mathrm{res}} E^1_{*,*,0} \cong \mathrm{dmf}_{*,*} \oplus \mathrm{dmf}_{*-1,*,-1}$$

sends 1 to 1 and $\tau$ to $\tau + h_1$ analogously to what we did in (3.12). We obtain a spectral sequence converging to $\mathrm{djF}_{**\infty}$, which we will proceed to compute.

We have the following integral result:

**Proposition 5.5.** *There is a multiplicative spectral sequence, called the algebraic Atiyah-Hirzebruch spectral sequence (AAHSS), convering to* $\mathrm{djF}_{**\infty}$ *with*

$$E^1_{***} = \mathrm{dmf}_{**}[z, \tau]/(\tau^2 - h_1\tau, z\tau - h_1\tau),$$

*where* $\tau = \mathrm{tr}(1) \in E^1_{1,2,0}$. *Moreover,* $d_1(\tau) = 0$ *and* $d_1(z) = \tau + h_1$.

*Proof.* By the multiplicativity of the exact couple, we know that

$$\tau^2 = \mathrm{tr}(1)^2 = \mathrm{tr}(1)\,\mathrm{res}(\tau).$$

Whatever $\mathrm{res}(\tau)$ is, it will detect $\eta$ by the compatibility with the descent spectral spectral sequences and (3.12), and so it must be nontrivial. However, there is only one nontrivial element in $\mathrm{dmf}_{1,1}$, and that is $h_1$.

The differential follows directly from the exact couple. □

Finally, also the modified exact couple (3.15) has an algebraic version converging to $\overline{\mathrm{djF}}_{**\infty} = \mathrm{dmf}_{**} \times_{\mathrm{djF}_{**0}} \mathrm{djF}_{**\infty}$, and its spectral sequence has the signature

$$E^1_{***} = \mathrm{dmf}_{**}[x, y]/(x^3 - y^2) \Longrightarrow \overline{\mathrm{djF}}_{**\infty}.$$

Its $D^1$-terms are $\overline{\mathrm{djF}}_{**m}$ and therefore themselves the $E^2$-terms of descent spectral sequences converging to $\mathrm{tmf}_*(P^m)$; moreover, any differential $d^n$ in it converge to differentials in the Atiyah-Hirzebruch spectral sequence for $\mathrm{tmf}_*(P^\infty)$.

Since our aim is to compute $\mathrm{djF}_{**,m}$ for all $m$, and not just $\mathrm{djF}_{**\infty}$, we need to modify the spectral sequences above by truncating the filtration (5.3) at $m$ (for all $m$) and introducing an additional grading (yes, a fourth one!). This is very reminiscent of the mechanics behind synthetic spectra. We will concentrate on the case at hand, although the method would work for any multiplicative unrolled exact couple.

**Theorem 5.6.** *There exists a multiplicative spectral sequence, called the unstable algebraic Atiyah-Hirzebruch spectral sequence (UAAHSS)*

$$E^{1,m}_{s,t,u} = \mathrm{dmf}_{**}[\tau, z, a]/(\tau^2 - \tau h_1, z\tau - zh_1) \Longrightarrow \mathrm{djF}_{2t-s,s,m}$$

*It satisfies the following properties:*

(1) *There is a map of spectral sequences s from the UAAHSS to the AAHSS of Prop. 5.5 collapsing the m-grading. This map is surjective on $E_1$ with kernel $(a - 1)$.*
(2) *If $x \in E_1$ of the UAAHSS is such that $s(x)$ is a $(k-1)$-cycle and $d_k(s(x)) = y$ then $x$ is a $(k-1)$-cycle in the UAAHSS and $d_k(x) = a^k y$.*

*In particular, the $d^1$ differential is completely determined by $d^1(z) = a(\tau + h_1)$.*

*Proof.* Let $D^1_t = \mathrm{djF}_{n,s,t}$ and $E^1_{***}$ as in the AAHSS. Then $(D^1_*, E^1_*)$ is the associated multiplicative exact couple. Denote by $(D^{1,m}, E^{1,m})$ the truncated exact couple with

$$D^{1,m}_{s,t,u} = D^1_{s,t,\min(m,u)} \quad \text{and} \quad E^{1,m}_{s,t,u} = \begin{cases} E^1_{s,t,u}; & u \leq m \\ 0; & u > m. \end{cases}$$

Clearly, the associated spectral sequence converges strongly to $D^{1,m}_{s,t,\infty} = \mathrm{djF}_{2t-s,s,m}$, but it is no longer multiplicative. However, the multiplicative structure of $(D^1, E^1)$ induces a multiplicative structure on
$$\left(\bigoplus_m D^{1,m}_{***}, \bigoplus_m E^{1,m}_{***}\right)$$
with $D^{1,m}_{s,t,u} \otimes D^{1,m'}_{s',t',u'} \to D^{1,m+m'}_{s+s',t+t',u+u'}$ and similarly for $E$. The associated spectral sequence is the UAAHSS and converges to $\bigoplus_m \mathrm{djF}_{**,m}$.

The element $\tau \in E^1_{1,2,0}$ lifts to an element in $E^{1,0}_{1,2,0}$, and the element $z \in E^1_{0,2,1}$ lifts to an element in $E^{1,1}_{0,2,1}$. We also have the element $a \in E^{1,1}_{0,0,0}$ which maps to 1 under the isomorphism $E^{1,1}_{0,0,0} \cong \mathrm{dmf}_{0,0}$. For $m > 0$, the term $E^{1,m}$ is freely generated over $\mathrm{dmf}_{**}$ by the monomials of degree $m$ in $a$ and $z$, while for $m = 0$, $E^{1,0}$ is just $\mathrm{djF}_{*,*,0}$. The description of the $E^1$-term of the UAAHSS follows.

Consider the map $s\colon D^{1,*}_{s,t,u} = \bigoplus_{m\geq 0} D^1_{s,t,\min(m,u)} \to D^1_{s,t,u}$ given on every summand by multiplication by $a^{u-\min(m,u)}$. This yields a map of exact couples from the UAAHSS to the AAHSS
$$(D^{1,*}_{s,t,u}, E^{1,*}_{s,t,u}) \to (D^1_{s,t,u}, E^1_{s,t,u}),$$
surjective on $E^1$ and mapping $a$ to 1. Thus it induces an isomorphism on the quotient by $(a-1)$ and hence an isomorphism of exact couples, and the resulting spectral sequences are isomorphic.

The claim about differentials follows from this and quadegree considerations. □

Note that $E^{1,m}_{s,t,u} = 0$ for $u > m$, and the only generator that has $u \neq m$ is the stabilization class $a$ with $m = 1$, $u = 0$. The degrees $u$ and $m$ of a class in the $E^1$-term of the UAAHSS should be interpreted as follows: every class $x$ not divisible by $a$ has $u = m$, which is the index (topologically, the cell) where it first appears ("is born"). Denote again by $s$ the map UAAHSS $\to$ AAHSS. If $s(x)$ is a nontrivial permanent cycle then all multiples $a^k x$ are nontrivial permanent cycles. If, on the other hand, $s(x)$ is a nontrivial $(k-1)$-cycle and the target of a $d^k$-differential then $x, \ldots, a^{k-1}x$ are still nontrivial infinite cycles representing derived Jacobi forms, but $a^k x$ dies. Thus the difference $m - u$ of any class $y$ denotes the "age" of that class.

We will now invert 2, which simplifies the spectral sequence significantly. We remind the reader (cf. [Bau08]) that
$$\mathrm{dmf}_{**}[\tfrac{1}{2}] \cong \mathbf{Z}[\alpha, \beta, c_4, c_6, \Delta]/(\alpha^2, c_4\alpha, c_4\beta, c_6\alpha, c_6\beta, 1728\Delta - c_6^2 + c_4^3)$$
as displayed in Figure 5.1. In this chart, a box represents a copy of $\mathbf{Z}[\tfrac{1}{2}]$ and a dot a copy of $\mathbf{Z}/3\mathbf{Z}$. A line of slope $\tfrac{1}{3}$ denotes multiplication by $\alpha$.

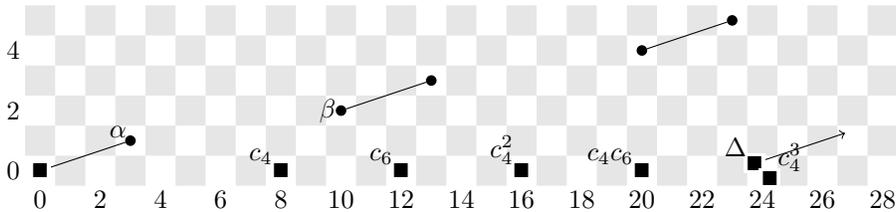

FIGURE 5.1. The bigraded ring $\mathrm{dmf}_{**}[\tfrac{1}{2}]$

**Lemma 5.7.** *The $E_2$-term of the UAAHSS of Thm. 5.6 after localization away from $2$ is given by*

$$E^{2,*}_{***} = \mathrm{dmf}_{**}[\tfrac{1}{2}, \tau, x, y, a]/(\tau^2, x\tau, y\tau, y^2 - x^3, a\tau)$$

*with $\tau$, $a$ as before and $x$, $y$ mapping to $z^2$ and $z^3$, respectively, in the $E^1$-term. The class $y$ is an infinite cycle, and there is a hidden extension $y\tau = 2ax\alpha$.*

*Proof.* The $E^2$-term follows from the differential $d^1(z) = a\tau$. (Recall that $h_1 = 0$ since $2h_1 = 0$.) It is displayed in Figure 5.2. We display it modulo the ideal $(c_4, c_6, \Delta)$ since it consists of permanent cycles.

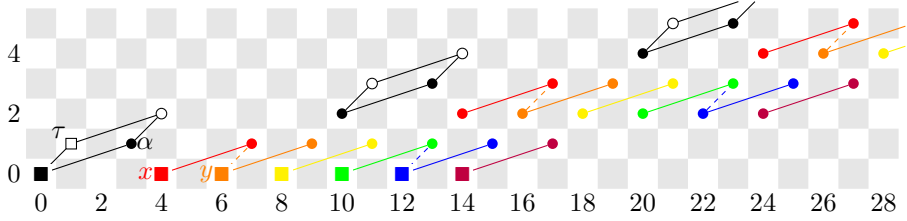

FIGURE 5.2. The $E^2$-term of the UAAHSS for $\mathrm{djF}_{***}[\tfrac{1}{2}]$

In that chart, the horizontal axis is the total dimension $n = t - s$, the vertical axis is the cohomological degree $s$, and RMA color code modulo 9 is used for the indices $u$ and $m$ (black = 0, brown = 1, red = 2, orange = 3, yellow = 4, green = 5, blue = 6, purple = 7, gray = 8). More precisely, a symbol with color code $i$ represents

- ■ $= \mathbf{Z}[\tfrac{1}{2}, a]$
- □ $= \mathbf{Z}[\tfrac{1}{2}]$
- ● $= \mathbf{Z}/3\mathbf{Z}[a]$
- ○ $= \mathbf{Z}/3\mathbf{Z}$,

with the generator 1 placed at bidegree $(u, m) = (i, i)$. A line of slope 1 represents a multiplication by $\tau$, and a line of slope $1/3$ one by $\alpha$.

The differentials $d^i$ have degree $(s, t, m, u) = (1, 0, 0, -i)$. Thus there are no possible targets for $d^i(y)$ for any $i \geq 2$. Thus there exists a lift $\tilde{y} \in \mathrm{djF}_{6,0,3}$ of $y$ such that $\mathrm{res}(\tilde{y}) = y$.

To see the exotic extension, note that $\mathrm{tr}(z)$ is a lift of $\tau$ to $\mathrm{djF}_{1,1,0}$. Thus $\tilde{y}\,\mathrm{tr}(z) = \mathrm{tr}(z\,\mathrm{res}(\tilde{y})) = \mathrm{tr}(z \cdot z^3) = \mathrm{tr}(z^4)$. We now use the interpretation of the topological spectral sequence of Prop. 3.14, or rather its $\overline{\mathrm{TMF}}$-variant, as the Atiyah-Hirzebruch spectral sequence for $P^m$. In particular, the topological transfer map $\mathrm{tr}\colon \Sigma^7 \mathrm{TMF} \to \mathrm{TJF}_3$ is induced by the attaching map $S^7 \to P^3$ of the 8-cell of $P^4$, which is the same as the attaching map of the 6-cell to the 2-cell of $\mathbf{C}P^3$, which is $2\nu$ with $\nu \in \pi_3(S^0)$ the Hopf map [Mos68]. (That it is nontrivial in mod-3 cohomology also follows from the unstable Steenrod algebra structure of $H^*(\mathbf{C}P^3, \mathbf{F}_3)$.)

We conclude that $\mathrm{tr}(z^4)$ has to be a class detecting $2x\nu$. Considering quadegrees and the fact that $x\nu$ is uniquely detected by $x\alpha$, we conclude that $\mathrm{tr}(z^4) = 2ax\alpha$. □

We proceed to compute the higher differentials in the UAAHSS. Since this spectral sequence arises from the filtration by powers of $a$ of the Hopf algebroid $H_{\mathcal{W},*}$ of (4.4) and in that Hopf algebroid,
$$\eta(x) = x - a^2 r,$$
and since furthermore $\alpha$ is represented by $[r]$ in the bar complex, we have
$$d_2(x) = a^2 \alpha.$$

By the Leibniz rule, we find that $d_2(x^i y^j) = \pm a^2 x^{i-1} y^j \alpha$ for $3 \nmid i$. The class $x\alpha$ survives to a class $\gamma$ of dimension 7 and index 1, supporting a multiplicative extension $\gamma\alpha = a^2\beta$ coming from the Massey product shuffle
$$\gamma\alpha = \langle a^2, \alpha, \alpha \rangle \alpha = a^2 \langle \alpha, \alpha, \alpha \rangle = a^2 \beta.$$

Since $a^2\alpha$ is a boundary, its $\gamma$-multiple $\gamma a^2 \alpha = a^4 \beta$ cannot survive either, and the only eligible differential to achieve this is
$$d_4(\{x^2 \alpha\}) = \pm a^4 \beta.$$

(A $d_3$ differential is impossible since all differentials preserve the $m$-grading, i.e. the color.)

Chart 5.3 shows the resulting $E_5$-term. The same conventions as in Chart 5.2 apply, with the addition that $k$ nested circles represent a copy of $\mathbf{Z}/3\mathbf{Z}[a]/(a^k)$. The classes $b_2$, $b_3$, and $b_4$ are chosen such that their restrictions are $3x$, $y$, and $3x^2$, respectively.

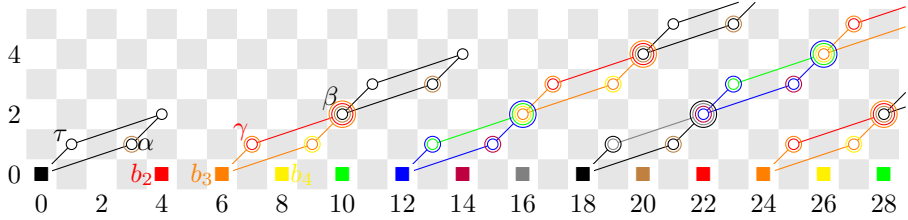

Figure 5.3. The $E_\infty$ term of the unstable algebraic Atiyah-Hirzebruch spectral sequence for $H^*(\mathcal{W}, L_{*,*})[\tfrac{1}{2}]$

The class $b_3$ is a permanent cycle, and no more differentials are possible for degree reasons.

*Proof of Thm. 4.6.* We read off from this chart that $(\mathrm{djF}_{***})[\tfrac{1}{2}]$ is generated by the classes $\alpha, \beta, c_4, c_6, \Delta$ that generate $(\mathrm{dmf}_{**})[\tfrac{1}{2}]$ as well as the classes $b_2$, $b_3$, $b_4$ in tridegrees $|b_i| = (2i, 0, i)$, the class $a$ in tridegree $(0, 0, 1)$, $\tau$ in tridegree $(1, 1, 0)$ as well as $\gamma$ in tridegree $(7, 1, 2)$.

We can choose the scaling of the classes $b_i$ so that they stabilize to the classes of the same name in $(\mathrm{djF}_{*,0,\infty})[\tfrac{1}{2}]$, cf. Thm. 4.5, so that the relations
$$a^4 c_4 = b_2^2 - 24 b_4, \quad a^6 c_6 = -b_2^3 + 36 b_2 b_4 - 216 b_3^2$$
and
$$a^{12} \Delta = \frac{1}{4} b_2^3 b_3^2 - 27 b_3^4 + 9 b_2 b_3^2 b_4 + \frac{1}{4} b_2^2 b_4^2 - 8 b_4^3$$
continue to hold. (The powers of $a$ are determined by the level of the equations, viz. $m = 2$, $m = 3$, and $m = 6$, respectively.)

The following relations hold:

- The relations involving $c_4$, $c_6$, and $\Delta$: $1728\Delta = c_6^2 - c_4^3$, $c_4\alpha = c_4\beta = c_4\gamma = c_6\alpha = c_6\beta = c_6\gamma = 0$, and the relations for $a^4c_4$, $a^6c_6$, and $a^{12}\Delta$ above;
- The relations present in $\mathrm{dmf}_{**}$: $\alpha^2 = 0$, $3\alpha = 3\beta = 0$;
- The relations coming from the kernel of the differentials on the $(s = 0)$-line: $b_2\alpha = b_4\alpha = b_2\beta = b_4\beta = 0$;
- The relations involving the new class $\gamma$: $3\gamma = a^2\gamma = \gamma^2 = b_2\gamma = b_4\gamma = 0$ and $\gamma\alpha = a^2\beta$;
- The $a$-order of classes in positive degree, all implied by $a^2\alpha = 0$;
- the relations involving $\tau$: $\tau^2 = \tau a\tau b_2 = \tau b_4 = 0$, $\tau b_3 = 2a\alpha$.

The claimed results follow. □

**Corollary 5.8.** *The associated graded of the filtration of $\mathrm{djF}_{*,*,\infty}$ is given by*

$$\mathrm{gr}\,\mathrm{djF}_{*,*,\infty}[\tfrac{1}{2}] = \mathrm{mf}_*[\tfrac{1}{2}, b_2, b_3, b_4]/I$$

$$I = (b_2^2 - 24b_4, -b_2^3 + 36b_2b_4 - 216b_3^2, \tfrac{1}{4}b_2^3b_3^2 - 27b_3^4 + 9b_2b_3^2b_4 + \tfrac{1}{4}b_2^2b_4^2 - 8b_4^3)$$

*Proof.* The associated graded is given by $\overline{\mathrm{djF}}_{*,*,*}/(a)$ and thus follows directly from Thm. 4.6. □

Moving on to the descent spectral sequence, we will determine all differentials from the tmf-module structure of $\mathrm{tjF}_\infty$ alone, not using any conjectural ring spectrum structure. We note that all differentials respect the index $m$ – i.e., the colors – and thus the only possible differentials are those induced by the fundamental differentials $d_5(\Delta) = \beta^2\alpha$ and $d_9(\Delta^2\alpha) = 2\beta^5$ implied by the tmf-module structure. This is displayed in Figure 5.4. Here, the black square in bidegree $(24,0)$ and its $b_3$-multiples indicate that the class $\Delta$ (black) died, but both $3\Delta$ and $a^2\Delta$ survive, thus it represents $(a^2, 3)\mathbf{Z}[\tfrac{1}{2}, a]$. Similarly, the black square in bidegree $(48, 0)$ is $(a^4, 3)[\tfrac{1}{2}, a]$. The classes $c_4, c_6$ are infinite cycles and not shown.

## 6. Computations at $p = 2$

We employ the same ideas to compute $\mathrm{djF}_{*,*,\infty}$ localized at the prime 2 as at the prime 3.

By completing the cube in $x$ in the polynomial $E(\underline{a}, x, y)$ (cf. [Bau08, Section 7]), the Hopf algebroid $H_{(2)} = (A_{(2)}, \Gamma_{(2)})$ is equivalent to $\tilde{H} = (\tilde{A}, \tilde{\Gamma})$ with

$$\tilde{A} = \mathbf{Z}_{(2)}[a_1, a_3, a_4, a_6], \quad \tilde{\Gamma} = \tilde{A}[s, t],$$

where $r = \tfrac{1}{3}(s^2 + a_1 s)$, and thus $H_{\mathcal{W},\infty} \simeq \tilde{H}_{\mathcal{W},\infty} = (\tilde{A}_\mathcal{W}, \tilde{A}_\mathcal{W}[s,t])$ with

$$\tilde{A}_\mathcal{W} = (\tilde{A}[x,y])/(E(a_1, 0, a_3, a_4, a_6, x, y)).$$

*Proof of Thm. 4.7.* We apply Theorem 5.1 for $A'_\mathcal{W} = \tilde{A}_\mathcal{W}/(y)$. We have that

$$g\colon \tilde{A}_\mathcal{W} \to A'_\mathcal{W} \otimes_{\tilde{A}_\mathcal{W}} \tilde{A}_\mathcal{W}[s,t] = A'_\mathcal{W}[s,t]$$

has $g(y) = -sx + sr - t = -sx + \tfrac{1}{3}(s^3 + a_1s^2) - t$. We can lift $g$ to a map

$$\tilde{g}\colon \tilde{A}[x,y] \to \tilde{A}[x,y,s,t]/(y)$$

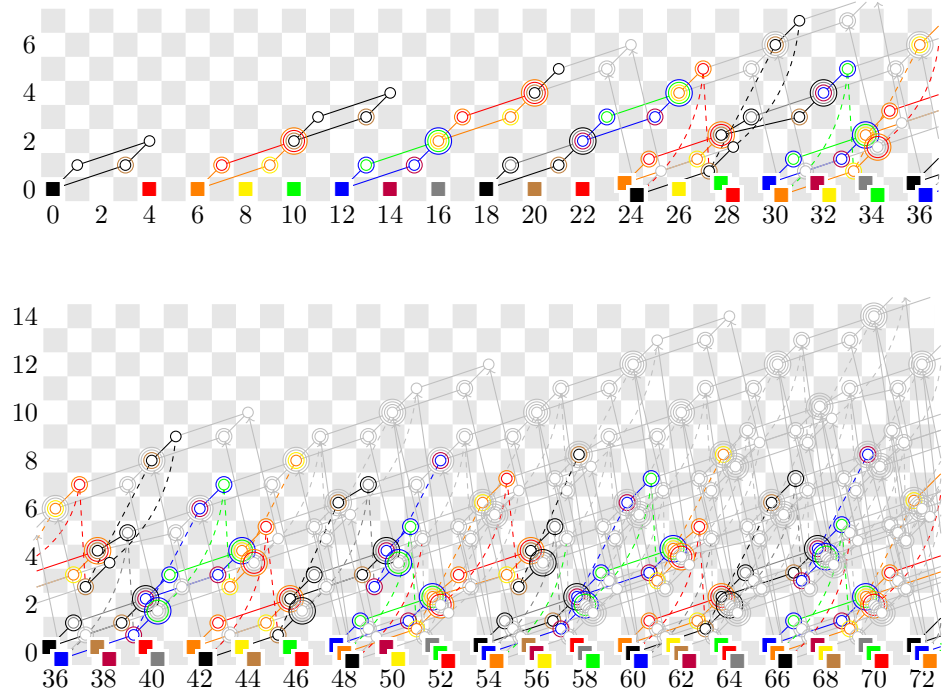

FIGURE 5.4. The descent spectral sequence converging to $\pi_* \operatorname{tjF}_*[\frac{1}{2}]$

which shows that the target is a free module over the source on the basis $\{s^i\}$, thus is faithfully flat. Hence so is the quotient by the (invariant) ideal generated by $E(a_1, 0, a_3, a_4, a_6, x, y)$. Thus $\tilde{H}_{\mathcal{W},\infty}$ is equivalent to $(A'_{\mathcal{W}}, \Gamma'_{\mathcal{W}})$ with

$$A'_{\mathcal{W}} = \mathbf{Z}_{(2)}[a_1, a_3, a_4, a_6, x]/(x^3 + a_4 x + a_6) \cong \mathbf{Z}_{(2)}[a_1, a_3, a_4, x]$$

and

$$\Gamma'_{\mathcal{W}} = A'_{\mathcal{W}}[s, t]/(sx - \frac{1}{3}s^2(s + a_1) + t) \cong A'_{\mathcal{W}}[s].$$

In $(A'_{\mathcal{W}}, \Gamma'_{\mathcal{W}})$, we have that

$\eta_R(a_1) = a_1 + 2s$

$\eta_R(a_3) = a_3 + a_1 r + 2t = a_3 + a_1 r + 2sr - 2sx = a_3 + \frac{1}{3}(a_1 + 2s)(s^2 + a_1 s) - 2sx$

$\eta_R(a_4) = a_4 - a_3 s - a_1(t + rs) - 2st + 3r^2 = a_4 - a_3 s + (a_1 + 2s)sx - 3r^2$

$\qquad = a_4 - a_3 s + (a_1 + 2s)sx - \frac{1}{3}(s^2 + a_1 s)^2$

$\eta_R(x) = x - r = x - \frac{1}{3}(s^2 + a_1 s).$

This shows that a further application of Thm. 5.1 for $A''_{\mathcal{W}} = A'_{\mathcal{W}}/(x)$ is possible since $A''_{\mathcal{W}} \otimes_{A'_{\mathcal{W}}} \Gamma'_{\mathcal{W}}$ is free of rank 2 over $B'$ with basis $\{1, s\}$. We thus obtain the equivalent Hopf algebroid $(A''_{\mathcal{W}}, \Gamma''_{\mathcal{W}})$ with

$$A''_{\mathcal{W}} = \mathbf{Z}_{(2)}[a_1, a_3, a_4]; \quad \Gamma''_{\mathcal{W}} = A''_{\mathcal{W}}[s]/(s^2 + a_1 s)$$

and
$$\eta_R(a_1) = a_1 + 2s$$
$$\eta_R(a_3) = a_3$$
$$\eta_R(a_4) = a_4 - a_3 s.$$

To compute the cohomology of this Hopf algebroid, we first divide out by the invariant ideal $I = (2, a_1)$. The resulting Hopf algebroid $(B_1, \Gamma_1)$ is given by
$$B_1 = \mathbf{F}_2[a_3, a_4]; \quad \Gamma_1 = B_1[s]/(s^2).$$
We have
$$\mathrm{Ext}_{(B_1, \Gamma_1)}(B_1, B_1) = \mathbf{F}_2[h_1, a_3, a_4^2]/(a_3 h_1),$$
where $h_1$ is in filtration 1 and total degree 1, represented in the bar construction by $[s]$.

Running the algebraic Bockstein spectral sequence to get the 2-local cohomology, the fundamental differential is given by $d^1(a_1) = 2s$, yielding
$$H^*(A''_\mathcal{W}, \Gamma''_\mathcal{W}) = \mathbf{Z}_{(2)}[h_1, a_1^2, a_3, a_3 a_1 + 2a_4, a_4^2 + a_4 a_3 a_1]/(2h_1, a_3 h_1)$$
or, in terms not referring to a surrounding chain complex,
$$\mathbf{Z}_{(2)}[h_1, b_2, b_3, b_4, b_8]/(2h_1, b_3 h_1, b_4 h_1, b_4^2 - b_2 b_3^2 - 4b_8).$$
with $c_4 = b_2^2 - 24 b_4$, $c_6 = -b_2^3 - 216 b_3^2 + 36 b_2 b_4$, $\Delta = -b_8 b_2^2 + 9 b_2 b_3^2 b_4 - 27 b_3^4 - 8 b_4^3$. □

In the descent spectral sequence for TJF, the differential $d_3(b_2) = h_1^3$ is forced by the relation $\eta^4 = 0$ in the sphere. This argument is valid in the (nonmultiplicative, but connective) version
$$\mathrm{tjF} \simeq \mathrm{tmf} \wedge P^\infty,$$
because $b_2$ is the only class that could support a differential, and thus it holds in TJF after inverting $\Delta$. By multiplicativity in TJF, this implies $d_3(b_2^{2n+1}) = b_2^{2n} h_1^3$ for all $n$. For dimensional reasons, $b_3$, $b_4$, and $b_8$ are infinite cycles and thus $d_3(b_8^k b_2^{2n+1}) = b_8^k b_2^{2n} h_1^3$ and we obtain
$$E_4 \cong \mathbf{Z}_{(2)}[h_1, 2b_2, b_3, b_2^2, b_4, b_2 b_3, b_2 b_4, b_8, \Delta^{-1}]/I$$
with
$$I = (2h_1, h_1^3, b_3 h_1, b_4 h_1, 4b_8 - b_2 b_3^2 + b_4^2).$$
For degree reasons, $E_4 \cong E_\infty$ and no multiplicative extensions are possible, and we conclude:

**Corollary 6.1.**
$$\pi_*(\mathrm{TJF}_\infty)_{(2)} \cong \mathbf{Z}_{(2)}[\eta, x_2, x_3, x_4, x_4', x_5, x_6, x_8, \Delta^{-1}]/I$$
*with*
$$I = (2\eta, \eta^3, x_2 \eta, x_3 \eta, x_4' \eta, x_5 \eta, x_6 \eta,$$
$$4x_4 - x_2^2, 2x_5 - x_2 x_3, 2x_6 - x_2 x_4', 2x_3 x_4 - x_2 x_5, 2x_4 x_4' - x_2 x_6,$$
$$x_5 x_6 - x_3 x_4 x_4', 4x_8 - (x_4')^2 + x_3 x_5).$$
*and*
$$\Delta = -x_8 x_4 + 9 x_5 x_3 x_4' - 27 x_3^4 - 8(x_4')^3.$$

This looks complicated; however, additively,

$$\pi_*(\mathrm{tjF}_\infty) \cong \bigoplus_{n\geq 0} \Sigma^{16n}\left(ko_* \oplus \bigoplus_{i=1}^{\infty} \Sigma^{6i}ku_*\right),$$

and this can be expressed even more concisely by saying that

$$\pi_*(\mathrm{tjF}_\infty) \cong ko_*(X)$$

where

$$X = \bigvee_{n\geq 0} \Sigma^{16n}\left(S^0 \vee \Sigma^6 C(\eta) \vee \Sigma^{12} C(\eta)\right).$$

The 16-periodic summands are given generated by $x_8^n$, the inner $ko_*$-summands by the subring generated by $x_2$ and $x_3'$, and the $ku$-factors by the subrings generated by $x_3^i, x_3^{i-1}x_4, x_3^{i-1}x_5, x_3^{i-1}x_6$.

We will not determine the full $\mathrm{tmf}_*$-algebra structure, but note that on the 0-line, the topological modular forms are mapped as $c_4 = b_2^2 - 24b_4 \mapsto x_4 - 24x_4'$ and $2c_6 = -2b_2^3 - 432b_3^2 + 72b_2b_4 \mapsto x_2x_4 - 432x_3^2 + 36x_2x_4'$.

## Appendix A. Equivariant homotopy theory and the transfer

In this appendix, we recall some equivariant homotopy theory and describe, in particular, different ways to obtain the transfer. To ensure consistency, we will take as only input equivariant Atiyah duality and construct all other maps and equivalences from there. All material here is known; only the viewpoint might be somewhat new.

Let $G$ be a compact Lie group. When speaking about $G$-spectra, we mean the $\infty$-category $\mathrm{Sp}_G$ of genuine $G$-spectra with respect to a full universe. We view pointed $G$-spaces implicitly as $G$-spectra via $\Sigma^\infty$. As a matter of convention, all our actions will be left actions in this appendix; if we speak of right actions, we act indeed by the inverse from the right.

The functor $\Sigma^\infty$ factors through the (naive) stabilization of $G$-spaces, called *naive $G$-spectra*. The corresponding right adjoint from genuine $G$-spectra $X$ to naive $G$-spectra records the $G$-spaces $X(V)$ for *trivial* $G$-representations $V$ (and the structure maps between them). The forgetful functor from genuine $G$-spectra to spectra with $G$-action factors in this way through naive $G$-spectra.

Let $G$ be a compact Lie group. Then $G_+$ is a $G \times G$-space, where the first $G$-factor acts by multiplication on the left and the second factor by multiplication with inverses on the right. We also denote by $G_+$ the corresponding $G \times G$-equivariant suspension spectrum $\Sigma^\infty_+ G$. The Spanier-Whitehead dual $DG_+ = F(G_+, S)$ is also a $G \times G$-equivariant spectrum.

Let $\pi_i \colon G \times G \to G$ be the projection maps and $\iota_i \colon G \to G \times G$ the inclusion maps ($i = 1, 2$). Denote by $G_i < G \times G$ the direct factor subgroups $\iota_i(G)$. Any pointed $G$-space or $G$-spectrum $X$ can be turned into a $G \times G$-spectrum $X^i$ in two ways, namely as $\pi_i^* X$. We consider the space $G$ as a $G \times G$-space with the left and right multiplication actions, and its $G \times G$-equivariant suspension spectrum $G_+ = \Sigma^\infty_+ G$.

We record the following straightforward coinduction result:

**Lemma A.1.** *Let $E$ be a $G \times G$-spectrum. There is a natural equivalence of naive $G$-spectra*

$$F(G_+, E)^{G_i} \simeq \Delta^* E \quad (i = 1, 2),$$

*where* $\Delta \colon G \to G \times G$ *is the diagonal.*

*Proof.* Let us first consider the case where $E$ is a $G \times G$-space. Then the evaluation map
$$F(G_+, E)^{G_1} \to \Delta^* E, \quad \phi \mapsto \phi(1),$$
is a $G$-equivalence: for the $G_2$-action $\phi \mapsto \phi.\gamma$ of $\gamma$ on $F(G_+, E)^{G_1}$, we have
$$(\phi.\gamma)(-) = (1,\gamma).\phi(- \cdot \gamma) \mapsto (\gamma,\gamma).\phi(1).$$
In the spectrum case, we have for a trivial $G$-representation $V$:
$$F(G_+, E)^{G_1}(V) = F(G_+, E(V))^{G_1} = \Delta^*(E(V)) = (\Delta^* E)(V). \qquad \square$$

Let $S^{\mathrm{Ad}\,G} = (\pi_1)^* S^{\mathrm{ad}\,G}$, where $S^{\mathrm{ad}\,G}$ is the $G$-spectrum which is the suspension spectrum of the one-point compactification of the adjoint representation of $G$.

The pivotal construction for both the Adams isomorphism and the transfer is the following special case of equivariant Atiyah duality:

**Lemma A.2.** *For any compact Lie group $G$, there is an equivalence of genuine $G \times G$-spectra*
$$G_+ \simeq S^{\mathrm{Ad}\,G} \wedge DG_+.$$

As a corollary, by taking fixed points with respect to the subgroup $G_2$ and using Lemma A.1, we obtain an equivalence of $G$-spectra
$$(G_+)^{G_2} \simeq (S^{\mathrm{Ad}\,G} \wedge DG_+)^{G_2} \cong F(G_+, S^{\mathrm{Ad}\,G} G)^{G_2} \cong S^{\mathrm{ad}\,G}.$$

*Proof.* By the Mostow-Palais (equivariant Whitney embedding) theorem, $G$ can be $G \times G$-equivariantly embedded into a $G \times G$-representation $V$. By the equivariant tubular neighborhood theorem [Bre72], there exists a $G \times G$-invariant tubular neighborhood $U$ of $G$ in $V$ together with a $G \times G$-diffeomorphism $\nu_{G \hookrightarrow V} \to U$, where $\nu_{G \hookrightarrow V}$ denotes the normal bundle of $G$ in $V$. The ($G \times G$-equivariant) Pontryagin-Thom collapse maps
$$S^V \to U \cup \{\infty\} \cong G^{\nu_{G \hookrightarrow V}}$$
assemble to a map of $G$-spectra
$$d_G \colon S \to G^{-TG}.$$
Composing with the Thom diagonal $G^{-TG} \to G^{-TG} \wedge G_+$ and taking duals yields a $G \times G$-map
$$(A.3) \qquad\qquad t_G \colon DG_+ \to G^{-TG},$$
which is an equivalence by [May96, Ch. XVI.8]. Note that the bundle $TG$ on $G$ is trivial as a left or right $G$-bundle by the (left or right) invariant framing. More precisely, the map
$$G \times T_e G \to TG, \quad (g, x) \mapsto (g, xg)$$
is a trivialization of $TG$, where $xg \in T_g G$ denotes the tangent vector obtained as the image of $x \in T_e G$ by the action of the derivative of right multiplication with $g$. Considering $TG$ as a $G \times G$-space by left and right multiplication, we see that the induced right $G$-action on $T_e G$ is trivial, while the left $G$-action is the adjoint representation (i.e. conjugation) on $T_e G$. Thus we have, $G \times G$-equivariantly, that $G^{-TG} \cong G_+ \wedge S^{-\mathrm{Ad}\,G}$. Combining with the equivalence (A.3) and smashing both sides with $S^{\mathrm{Ad}\,G}$ yields the claim. $\qquad \square$

From the previous lemma, we obtain for every $G$-spectrum $E$ a chain of equivalences

$$(A.4) \quad \begin{aligned}(E^1 \wedge G_+)^{G_1} &\simeq (E \wedge S^{\operatorname{Ad} G} \wedge DG_+)^{G_1} \simeq F(G_+, E \wedge S^{\operatorname{Ad} G})^{G_1} \\ &\simeq E \wedge S^{\operatorname{ad} G} \simeq (E^1 \wedge S^{\operatorname{Ad} G} \wedge G_+)_{hG_1}\end{aligned}$$

of $G$-spectra.

The Adams isomorphism is a generalizing this from $G$ to arbitrary complexes $X$ with free $G$-action. Below, we need the compatibility of equivariant Atiyah duality with the Adams isomorphism. For this reason (and to be more self-contained), we give a quick derivation of a special case of the Adams isomorphism ([LMSM86, SII.7]) from Lemma A.2.

**Proposition A.5** (Adams isomorphism). *If $X$ is a $G$-CW complex with a free $G$-action, then*

$$(A.6) \quad (E \wedge X_+)^G \simeq (E \wedge S^{\operatorname{Ad}(G)} \wedge X_+)_{hG}$$

*for every (genuine) $G$-spectrum $E$. This equivalence is functorial in $X$ and $E$. In the case of $X = G$, the equivalence is homotopic to (A.4).*

Note that if $E$ carries a trivial underlying $G$-action, the right-hand side simplifies to $\operatorname{res}_e^G E \wedge (S^{\operatorname{ad} G} \wedge X_+)/G$.

*Proof.* We will consider $G \times G$-equivariant spectra with the conventions discussed above Lemma A.1. For any pointed $G$-CW-complex $X$, the shear map

$$\operatorname{sh}\colon G \times X \to G \times X, \qquad (g, x) \mapsto (g, g^{-1}x)$$

induces a $G \times G$-isomorphism $s\colon G_+ \wedge X^1 \cong G_+ \wedge X^2$ (where $G \times G$ acts on both factors).

Now consider the $G \times G$-spectrum

$$Z = E^2 \wedge G_+ \wedge X^1_+.$$

Note that this is a free $G_1$-spectrum, with $G_1$-homotopy orbits

$$Z_{hG_1} \xrightarrow{\operatorname{id}_E \wedge s} (E^2 \wedge G_+ \wedge X^2_+)_{hG_1} \cong E \wedge X_+.$$

Taking $G$-fixed points on this gives one side of the Adams isomorphism, $(E \wedge X_+)^G$. On the other hand, taking $G_2$-fixed points of the same spectrum and using Lemmas A.2 and A.1

$$Z^{G_2} \cong (E^2 \wedge (S^{\operatorname{Ad} G} \wedge DG_+) \wedge X^1_+)^{G_2} \cong F(G_+, E^2 \wedge S^{\operatorname{Ad} G} \wedge X^1_+)^{G_2} \sim E \wedge S^{\operatorname{ad} G} \wedge X_+,$$

where the last equivalence is of naive $G$-spectra. Taking $G$-homotopy orbits of this yields the right hand side

$$(E \wedge S^{\operatorname{ad} G} \wedge X_+)_{hG}.$$

The desired map (A.6) is the exchange map $(Z^{G_2})_{hG} \to (Z_{hG_1})^G$ for homotopy orbits and fixed points. Fixed points commute with all homotopy colimits (see e.g. [MSZ23, Section 2.1]) and in particular with homotopy orbits. Thus, the exchange map is an equivalence.

The compatibility of the Adams isomorphism with Equation (A.4) is a somewhat lengthy but straightforward check.[3]　□

---

[3]For readers wanting to avoid this check: The space of functorial equivalences $(E \wedge G_+)^G \simeq E \wedge S^{\operatorname{Ad}(G)}$ is equivalent to that of self-equivalences of $S^{\operatorname{Ad}(G)}$. In particular, (A.4) and (A.6)

*Remark* A.7. There is an alternative proof of the Adams isomorphism in the form of Proposition A.5 using some $\infty$-categorical technology: The $\infty$-category of $G$-CW complexes with free $G$-action is by Elmendorf's theorem equivalent to the $\infty$-category of functors from the full subcategory $\{G\}$ of the orbit category $\mathrm{Orb}_G$ on $G$ to spaces. This is free generated under colimits from $\{G\}$. Thus, two colimit preserving functors of free $G$-spaces are equivalent if they are so on $\{G\}$. Thus, the equivalence (A.6) is induced by (A.4) since both sides preserve colimits in $Y$.

We use the following definition of the transfer, where we denote the underlying spectrum of a $G$-spectrum $E$ by $E^e$.

**Definition.** Given a $G$-spectrum $E$, we define the *transfer map* $E^e \wedge S^{\mathrm{ad}\, G} \to E^G$ as the composition
$$E^e \wedge S^{\mathrm{ad}\, G} \simeq F(G_+ \wedge S^{-\mathrm{ad}\, G}, E)^G \simeq F(DG_+, E)^G \to F(DS^0, E)^G \simeq E^G.$$
Here, we applied Lemma A.2 in the second step, and the third map is induced by the dual of the map $G_+ \to S^0$.

**Lemma A.8.** *For a $G$-spectrum $E$, the transfer map coincides with the composite*
$$E^e \wedge S^{\mathrm{ad}\, G} \simeq (E \wedge G_+)^G \to E^G,$$
*where the first equivalence is the Adams isomorphism, and the last map is induced by the map $G_+ \to S^0$.*

*Proof.* This follows from the compatibility in the last sentence of Proposition A.5. $\square$

The only relevant case for us is $G = T$, where $\mathrm{ad}\, G$ is the trivial 1-dimensional representation and the transfer above takes the form $\Sigma E^e \to E^T$.

We also give two definitions of a transfer map $\mathbf{C}P_+^m \to S^{-1}$. To that purpose, recall that $\rho$ denotes the tautological $T$-representation on $\mathbf{C}$.

(1) There is a cofiber sequence
$$T_+ \to S^0 \to S^\rho \to \Sigma T_+.$$
Desuspending by $\rho$ gives a map $\gamma \colon S^0 \to \Sigma^{1-\rho} T_+ \simeq \Sigma^{-1} T_+$. Taking $T$-homotopy orbits yields the map $\mathbf{C}P_+^\infty \simeq S_{hT}^0 \to (\Sigma^{-1} T_+)_{hT} \simeq S^{-1}$. We obtain the map $\mathbf{C}P_+^m \to S^{-1}$ by precomposing with $\mathbf{C}P_+^m \to \mathbf{C}P_+^\infty$.

(2) Using the Adams isomorphism, we obtain the map
$$(A.9) \qquad \mathbf{C}P_+^m \simeq \Sigma^{-1}(S((m+1)\rho)_+)^T \to \Sigma^{-1}(S^0)^T \to \Sigma^{-1} S^0 = S^{-1}.$$
Here, the last map is restriction or, equivalently, the desuspension of the composite
$$\iota \colon (S^0)^T \simeq (DS^0)^T \to (DT_+)^T \simeq S^0.$$

*Remark* A.10. The map $S^\rho \to \Sigma T_+$ in the first item agrees with the Pontryagin collapse map for the standard embedding $T \hookrightarrow S^\rho$. Thus, its $\rho$-th desuspension $\gamma \colon S^0 \to \Sigma^{-1} T_+$ agrees with the map $S^0 = DS^0 \to DT_+ \simeq \Sigma^{-1} T_+$. Indeed, the Atiyah duality equivalence from Lemma A.2 is induced from the Pontryagin collapse map.

---

can differ at most by a sign. Moreover, in Remark A.7, we give an alternative proof where the equivalence of (A.4) and (A.6) is already built in.

**Proposition A.11.** *These two definitions agree (up to sign). The fiber of this transfer agrees with the stunted projective space $\mathbf{C}P_{-1}^m$, i.e. the Thom spectrum of the negative of the tautological line bundle on $\mathbf{C}P^{m+1}$.*

*Proof.* Given the first definition of the transfer, the fiber of $\mathbf{C}P_+^\infty \to S^{-1}$ is $S_{hT}^{-\rho}$, which agrees with the stunted projective space $\mathbf{C}P_{-1}^\infty$. The fiber of $\mathbf{C}P_+^m \to S^{-1}$ is thus the fiber product $\mathbf{C}P_{-1}^\infty \times_{\mathbf{C}P_+^\infty} \mathbf{C}P_+^m$, which agrees with $\mathbf{C}P_{-1}^m$.

It remains to show the equality of transfers. Since the map $S(m\rho)_+ \to S^0$ factors as $S(m\rho)_+ \to S(\infty\rho)_+ \to S^0$ and the Adams isomorphism is functorial, the transfer $\mathbf{C}P_+^m \to S^{-1}$ in the second sense arises as the composition $\mathbf{C}P_+^m \to \mathbf{C}P_+^\infty \to S^{-1}$, with the second map again the transfer in the second sense. Thus, it suffices to prove the equality of transfers for $m = \infty$.

Consider the diagram

$$\begin{array}{ccccc}
\mathbf{C}P_+^\infty \simeq S(\infty\rho)_+/T \simeq S_{hT}^0 & \longrightarrow & \Sigma^{-1}T_+ \wedge S(\infty\rho)_+/T & \longrightarrow & \Sigma^{-1}T_+/T \simeq S^{-1} \\
\downarrow\simeq & & \downarrow\simeq & & \downarrow\simeq \\
\Sigma^{-1}S(\infty\rho)_+^T & \longrightarrow & \Sigma^{-2}(T_+ \wedge S(\infty\rho)_+)^T & \longrightarrow & \Sigma^{-2}(T)^T,
\end{array}$$

where the first horizontal map is induced by the map $\gamma\colon S^0 \to \Sigma^{-1}T_+$ discussed above, the second by the collapsing map $S(\infty\rho)_+ \to S^0$, and all vertical arrows are Adams isomorphisms. The composite in the top row is equivalent to the first definition of the transfer. In the second, we may alternatively first collapse $S(\infty\rho)$ and then apply $\gamma$. Thus, to identify the composite in the second row with the second description of the transfer, it remains to show that the map

$$\iota\colon (S^0)^T \simeq (DS^0)^T \to (DT_+)^T \simeq S^0$$

above agrees with the composite

$$(S^0)^T \xrightarrow{\gamma} (\Sigma^{-1}T_+)^T \simeq S^0.$$

Here, the last equivalence is the Adams isomorphisms or, equivalently, (A.4) for $E = S^0$. Thus, the agreement of the two maps follows from Remark A.10. □

As a last point, we want to prove a compatibility result about the transfer in the setting of *global spectra*. Roughly speaking, a global spectrum consists of a family of genuine $G$-spectra for all compact Lie groups $G$; see [LNP25] for a definition of the $\infty$-category $\mathrm{Sp}_{\mathrm{gl}}$ of global spectra in this spirit and [Sch18, Boh14] for (equivalent) other approaches. The important example for us is TMF, for which the structure of a global spectrum was constructed in [GLP24]. For a global spectrum $X$, we denote by $X_G$ the corresponding $G$-spectrum and by $X^G$ its fixed points.

**Lemma A.12.** *Let $X$ be a global spectrum. Then*

$$X_e \wedge \mathbf{C}P^m \simeq \Sigma^{-1}(X_T \wedge S((m+1)\rho)_+)^T \to \Sigma^{-1}X^T \xrightarrow{\mathrm{res}} \Sigma^{-1}X_e$$

*is equivalent to $X_e \wedge \mathrm{tr}_m$, where $\mathrm{tr}_m$ is the transfer as in (A.9). (Here, the first equivalence uses the Adams isomorphism.) Moreover, this equivalence is natural in $m$ and $X$.*

*Proof.* There is a unique colimit-preserving functor $L\colon \mathrm{Sp} \to \mathrm{Sp}_{\mathrm{glo}}$ sending the sphere to the sphere; this is left-adjoint to the restriction functor $Z \mapsto Z_e$. This is spelled out in model-categorical terms in [Sch18, Theorem 4.5.1], in whose proof it is also shown that the adjunction unit $Z \mapsto (LZ)_e$ is an equivalence.

The counit $L(X_e) \to X$ induces a diagram

$$\begin{array}{cccccc}
X_e \wedge \mathbf{C}P^m & \longrightarrow & \Sigma^{-1}X_e \wedge (S((m+1)\rho)_+)^T & \longrightarrow & \Sigma^{-1}X_e \wedge (S^0)^T & \longrightarrow & \Sigma^{-1}X_e \\
\downarrow \simeq & & \downarrow \simeq & & \downarrow \simeq & & \downarrow \simeq \\
(LX_e)_e \wedge \mathbf{C}P^m & \stackrel{\simeq}{\to} & \Sigma^{-1}((LX_e)_T \wedge S((m+1)\rho)_+)^T & \longrightarrow & \Sigma^{-1}(LX_e)^T & \stackrel{\mathrm{res}}{\Rightarrow} & \Sigma^{-1}(LX_e)_e \\
\downarrow \simeq & & \downarrow & & \downarrow & & \downarrow \simeq \\
X_e \wedge \mathbf{C}P^m & \stackrel{\simeq}{\longrightarrow} & \Sigma^{-1}(X_T \wedge S((m+1)\rho)_+)^T & \longrightarrow & \Sigma^{-1}X^T & \stackrel{\mathrm{res}}{\longrightarrow} & \Sigma^{-1}X_e.
\end{array}$$

The upper vertical arrows are special cases of the (natural) natural transformation $X \wedge F(-) \to F(X \wedge -)$ existing for every exact functor $F$ on $T$-spectra and every spectrum $X$, coming from the colimit-interchange map. In fact, these are equivalences in this case since $T$-fixed points commute with homotopy colimits. $\square$